\documentclass[10 pt, journal]{IEEEtran}

\IEEEoverridecommandlockouts                              

\usepackage{graphics} 
\usepackage{epsfig} 
\usepackage{times} 
\usepackage{amsmath} 
\usepackage{amssymb}  
\usepackage{amscd}
\usepackage{ifthen}
\usepackage{cite}
\usepackage{lscape}  
\usepackage[justification=raggedright]{caption}	
\usepackage{algorithm}
\usepackage{algpseudocode}
\usepackage{cite}

\def\scr#1{{\cal #1}}
\newcommand{\R}{\mathbb{R}}

\newtheorem{theorem}{Theorem}
\newtheorem{definition}{Definition}
\newtheorem{lemma}{Lemma}
\newtheorem{remark}{Remark}

\newtheorem{assumption}{Assumption}

\usepackage{color}
\usepackage{xcolor}
\definecolor{purple}{rgb}{0,0,0}
\definecolor{blue}{rgb}{0, 0, 0}

\begin{document}

\title{{\color{black}A Game-Theoretic Framework for Multi-Period-Multi-Company Demand Response
 Management in the Smart Grid
} \thanks{K. Alshehri is with the Systems Engineering Department, King Fahd University of Petroleum and Minerals (\texttt{kalshehri@kfupm.edu.sa}). J. Liu is with Department of Electrical and Computer Engineering, Stony Brook University
(\texttt{ji.liu@stonybrook.edu}). X. Chen is with the Department of Electrical, Computer, and Energy Engineering, University of Colorado Boulder  (\texttt{xudong.chen@colorado.edu}). T. Ba\c{s}ar is with the Department of Electrical and Computer Engineering and the Coordinated Science Laboratory, University of Illinois at Urbana-Champaign (\texttt{basar1@illinois.edu}). }
}


\author{Khaled Alshehri, ~Ji Liu, ~Xudong Chen, ~Tamer Ba\c sar\\
{\it{\color{red}To appear in IEEE Transactions on Control Systems Technology}}
}


\maketitle
\thispagestyle{empty}
\pagestyle{empty}

\begin{abstract}

By utilizing tools from game theory, {\color{black} we develop a novel multi-period-multi-company demand response framework considering the interactions between companies (sellers of energy) and their consumers (buyers of energy).} We model the interactions in terms of a Stackelberg game, where companies set their prices and consumers respond by choosing their demands. We show that the underlying game has a unique equilibrium at which the companies  maximize their revenues while the consumers maximize their utilities subject to their local constraints. Closed-form expressions are provided for the optimal strategies of all players. Based on these solutions, a power allocation game has been formulated, which is shown to admit a unique pure-strategy Nash equilibrium, for which closed-form expressions are also provided. {\color{black} This equilibrium is found under the assumption that companies can freely allocate their power across the time horizon, but we also demonstrate that it is possible to relax this assumption.} {\color{blue} We further provide a fast distributed algorithm for the computation of all optimal strategies using only local information.} {\color{black} We also study the effect of variations in the number of periods (subdivisions of the time horizon) and the number of consumers.} As a consequence, we are able to find an appropriate company-to-consumer ratio when the number of consumers participating in demand response {\color{black} exceeds some threshold}. Furthermore, we show, both analytically and numerically, that the multi-period scheme provides incentives for energy consumers to participate in demand response, {\color{black} compared to the single-period framework studied in the literature \cite{sabita}}. {\color{black}In our framework, we provide a condition for the minimum budgets consumers need, and carry out case studies using real life data to demonstrate the benefits of the approach, which show potential savings of up to $30\%$ and {\color{black} equilibrium prices that have low volatility.}}
 
\end{abstract}

\section{Introduction}

{\color{black}
One critical aspect of demand-side management (DSM) in the smart grid is demand response, which is defined as the response of consumers' demands to price signals from the utility companies (see \cite{DR,DSM,DSview} for tutorial discussions). Demand response allows companies to manage the consumers' demands, either directly (through direct load control) or indirectly (through pricing mechanisms). Demand response comes with great benefits, including -but not limited to- improving the electricity market efficiency \cite{DReff}. It also comes with challenges, particularly in its deployment \cite{challengesDR}. For an overview of the methodologies and the challenges of load/price forecasting and managing demand response in the smart grid, see \cite{DRprice}. A comprehensive survey on the pricing methods and optimization algorithms for demand response programs can be found in \cite{surveyDRM}. An overview of integrated demand response, where consumers participate in multiple energy systems is provided in \cite{reviewAE}.



 Using the framework of game theory, load adaptive pricing has been introduced decades ago  \cite{loadadapt}. {\color{blue}In this paper, we use game theory to design a multi-period-multi-company demand response management program at which companies and their consumers reach a unique equilibrium. At the equilibrium point, prices and demands are optimally chosen such that companies maximize their  revenues and consumers maximize their utility functions.} {\color{black} For the purpose of this paper, one can think of ``company" as a utility company serving households, businesses, and industrial consumers.} It is of critical interest to capture competition between companies, and hence we utilize the framework of, and tools from noncooperative game theory. We remark that such tools can be useful for the smart grid in various contexts \cite{survey}.


A considerable number of contributions have used game theory to analyze what happens in a smart grid where there are multiple sellers/utilities/retailers \cite{walidPHEV,trading,walidPHEV2,gao,twolevel,tansu,yaagoubi,tushar2,sabita,sabita2,sabita3,han} serving the same set of consumers. For example, analysis of how plug-in hybrid electric vehicles can sell back to the grid has been explored in \cite{walidPHEV,trading,walidPHEV2}. A similar analysis has also been carried out for electric bicycles \cite{gao}. A two-level game has been proposed in \cite{twolevel}. The authors in \cite{tansu} introduce a Stackelberg game to capture the interactions between electricity generator owners and a demand response aggregator. In \cite{yaagoubi}, a distributed game between energy consumers of different types has been designed while emphasizing individual preferences. Furthermore, in \cite{tushar2}, analysis of three-party energy management scheme between residential users, a shared facility controller, and the main power grid, has been conducted via a Stackelberg game. Among the contributions in the literature the ones most relevant to this paper are \cite{sabita} and \cite{sabita2}. A single-period Stackelberg game for demand response management with multiple utility companies has been proposed in \cite{sabita}, where consumers choose their optimal demands in response to prices announced by different utility companies. In \cite{sabita2}, an extension to the large population regime was carried out. Variations of \cite{sabita} to user-centric approaches were discussed in \cite{sabita3,han}. These works \cite{walidPHEV,trading,walidPHEV2,gao,twolevel,tansu,yaagoubi,tushar2,sabita,sabita2,sabita3,han} have demonstrated the usefulness and the power of game theory in capturing the interplay between buyers and sellers in the smart grid, but they are limited to single period setups.

In the smart grid, temporal variations play a critical role on both the supply side and the demand side. There are several papers in the literature that have addressed inter-temporal considerations in DSM and demand response \cite{amir,hazem,roh,zhudiff,PAR,collins,repeated,fourstage,dayahead,wei}, such as scheduling of appliances and/or storage \cite{amir,hazem,roh,zhudiff},  peak-to-average ratio reduction \cite{PAR,collins,repeated}, procurement issues \cite{fourstage}, and wholesale market price fluctuations \cite{dayahead,wei}. While the contributions in \cite{amir,hazem,roh,zhudiff,PAR,collins,repeated,fourstage,dayahead,wei} are important and reveal the importance of game theory for multi-period considerations in demand-side management, they are all limited to a single seller/utility/retailer case. 

}

{  {\color{black}The vast majority of demand response contributions are either limited to a single seller case, or a single period one. Furthermore, they primarily focus on either the utility-side or the consumer-side. Our goal here is to alleviate these limitations by developing a multi-period-multi-company demand response framework in which we address the interests and incentives for both utilities and their consumers in the smart grid. We achieve our goal by formulating and solving a Stackelberg game, which is a hierarchical game consisting of two kinds of players, leaders who act first, and these are utility companies in our framework, and followers who respond  to leaders' decisions, and these are price-responsive consumers. We prove that the proposed game admits a unique equilibrium at which companies find their revenue-maximizing prices and consumers choose their optimal demands that maximize their utility functions while taking into account their budget limitations and energy needs across the time horizon. We further propose a distributed algorithm to compute the equilibrium using only local information. The unique equilibrium is computed for the case in which the power available to sell for each company at each period is fixed. Nevertheless, by exploiting the closed-form solutions we derive, we are able to formulate a new power allocation game at which companies solve for allocations that further maximize their revenues, and also prove that it admits a unique equilibrium, and find its analytical expression. The equilibrium of the power allocation game reveals that companies find it optimal to sell the same amount of power at each period. This affirms that our game-theoretic framework aligns with the incentives of utility companies that prefer to minimize the Peak-to-Average ratio. Furthermore, we study what happens in the large population regime where the number of demand-responsive consumers becomes very large, and reveal that the number of companies needs to change appropriately, leading to an appropriate company-to-consumer ratio. {\color{blue} We also study what happens as the number of periods (subdivisions of the time horizon) grows, and show both theoretically and numerically, that consumers' utility increase as the number of periods increases, making multi-period demand response desirable for them.} Since we also address revenue-maximization for companies, this leads to a win-win situation. Furthermore, we provide a theoretical benchmark to measure whether or not consumers are spending more than what is necessary. We validate the applicability of our game to real life data. Numerical studies show that our benchmark leads to billing savings in excess of $10-30\%$, demonstrate the fast convergence of our distributed algorithm, and quantify the effect of the number of periods.  Our work captures the competition between companies, budget limitations at the consumer-level, and energy need for the entire time-horizon. \footnote{{\color{black}Some of the results in this paper were presented earlier in the conference paper \cite{mywork}, but this paper provides a much more comprehensive treatment of the work, such as the inclusion of the power allocation, asymptotic analysis, distributed algorithm, generalizations, and proofs.}} We stress that we make some simplifying assumptions to keep our analysis tractable, which makes it possible to reveal the main insights and gain deep understanding into the interplay between companies and their consumers. We also demonstrate that our framework has  desirable mathematical properties that make generalizations at both the consumers-level and companies-level possible, which we discuss in Section \ref{generalizations}. }

The remainder of the paper is organized as follows. Preliminaries from game theory are provided in Section \ref{prelim}. The problem is formulated in Section \ref{formulation}, and optimal prices and demands are obtained via a Stackelberg game in Section \ref{game1}. In Section \ref{game2}, a power allocation game at the companies side is formulated based on the closed-form solutions of the Stackelberg game. Next, we provide a distributed algorithm for the computation of all optimal strategies using local information in Section \ref{algorithm}. The asymptotic regimes are studied, in which the number of periods or the number of consumers grows in Section \ref{asymp}. Next, we present results on case studies using real demand response data in Section \ref{numerical}. Generalizations are discussed in Section \ref{generalizations}. Finally, we conclude the paper in Section \ref{conclusion} with a recap of its main points and identification of future directions. An appendix at the end provides details of proofs of the five theorems and some auxiliary results. 
\section{Preliminaries from Game Theory}\label{prelim}
A static $N$-person noncooperative game is comprised of players set, action sets, and utility functions. Let the players set be denoted by $\scr{N} := \{1,\dots,N\}$, where $N$ is the number of players. Each player has an action set $\scr{A}_i$, and the decision of player $i$ is denoted by ${\bf{a}}_i\in\scr{A}_i$. The vector of decisions taken by other players is denoted by $\mathbf{{\bf{a}}_{-i}}:=({\bf{a}}_1,\dots,{\bf{a}}_{i-1},{\bf{a}}_{i+1},\dots,{\bf{a}}_N)$. Each player $i$ aims to maximize his/her utility function $u_i({\bf{a}}_i,\mathbf{{\bf{a}}_{-i}})$. One key point is that the utility function of player $i$ depends not only on his/her actions, but also on the decisions made by other players. An equilibrium concept that is suitable for  such games is the Nash Equilibrium (NE), which is defined below. 

\begin{definition} The action vector $\mathbf{a^*} \in \scr{A}_1\times\dots\times\scr{A}_N$ constitutes a Nash equilibrium for the $N$-person static noncooperative game in pure-strategies if 
\begin{equation}u_i({\bf{a}}^*_i,\mathbf{a^*_{-i}}) \geq u_i({\bf{a}}_i,\mathbf{a^*_{-i}}) \hspace{0.2in} \forall {\bf{a}}_i\in\scr{A}_i, \ i \in {\cal N}. \label{NE}\end{equation} \end{definition}
Sometimes it would be beneficial to allow for hierarchy in the decision process. In such a case, there are two types of players, leaders and followers. The leaders' decisions are dominant, and the followers respond to the decisions taken by the leaders. This kind of {\color{blue}hierarchical} games is called Stackelberg games, and the corresponding solution concept is called the Stackelberg equilibrium. 
 The leaders have the privilege of choosing how to take their actions at the beginning of the game. However, they have to take into account how the followers would respond to these actions and how each leader's decision is influenced by the decisions of the other leaders. To be more precise, suppose that we have $K$ leaders and $N$ followers. Denote the followers set by $\scr{N} := \{1,\dots,N\}$, and the leaders set by $\scr{K} := \{1,\dots,K\}$, with action sets $(\scr{F}_i)_{i\in \scr{N}}$ and $(\scr{L}_j)_{j\in \scr{K}}$, respectively. Denote a generic action of leader $j$ by ${\bf{a}}_{j}\in\scr{L}_j$, and that of follower $i$ by  ${\bf{b}}_i\in\scr{F}_i$. The vector of actions taken by all leaders is denoted by $\mathbf{a}:=({\bf{a}}_1,\dots,{\bf{a}}_K)$. The utility of leader $j$ is denoted by $u_j({\bf{a}}_j,\mathbf{{\bf{a}}_{-j}},\mathbf{b(a)})$, where $\mathbf{{\bf{a}}_{-j}}$ denotes the decisions of the other leaders, and $\mathbf{b(a)}=({\bf{b}}_1(\mathbf{a}),\dots,{\bf{b}}_N(\mathbf{a}))\in \scr{F}_1 \times \dots \times \scr{F}_N$.

\begin{definition}The action vector $\mathbf{a^*} \in \scr{L}_1\times\dots\times\scr{L}_K$ is a Stackelberg Equilibrium strategy for all the $K$ leaders in pure-strategies if, for each $ j \in {\cal K}$,
\begin{eqnarray}u_j({\bf{a}}^*_j,\mathbf{a^*_{-j}},\mathbf{b^*(a^*)}) \geq u_j({\bf{a}}_j,\mathbf{a^*_{-j}},\mathbf{b^*}({\bf{a}}_j;{\bf{a^*_{-j}}}))\ \forall {\bf{a}}_j\in\scr{L}_j \label{SE}\end{eqnarray}
\end{definition}
where $\mathbf{b^*(a)} \in \scr{F}$ is the optimal response by all followers to the leaders' decisions, under the adopted equilibrium solution concept at the followers level. {\color{black} This solution concept is generally the Nash equilibrium, where followers play a Nash game. When there is no direct coupling between different followers, that is, other followers' decisions do not directly appear in the problem follower $i$ solves, they become independent, individual utility maximizers, which is the case we have in this work.} For a Stackelberg game, the pair ($\mathbf{a^*,b^*(a^*)}$) constitutes the equilibrium strategy. 

\section{Formulation of a Mathematical Model}
 \label{formulation}
\begin{figure}
\centering
\includegraphics[width=.9 \linewidth]{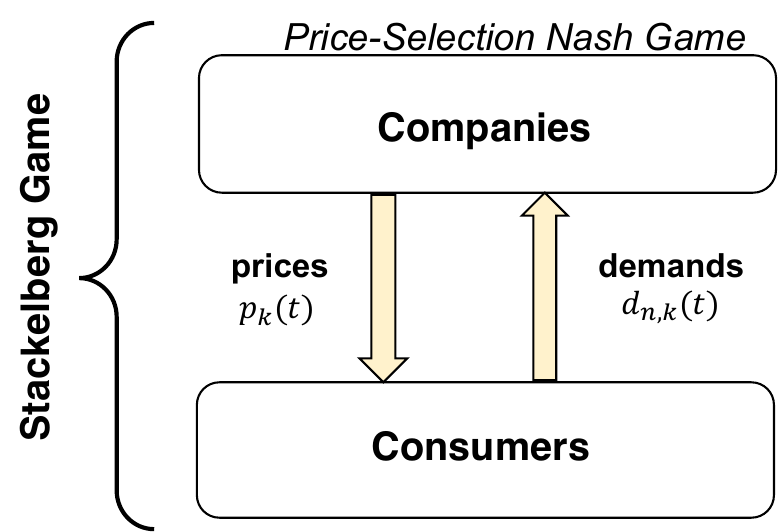}
\caption{The interaction between companies and their consumers. Companies play a price-selection Nash game. Then, consumers respond by choosing their demands independently of each other (the entire two-level interaction is a Stackelberg game).}
\label{sketch}
\end{figure}
Let $\scr{K}=\{1,2,\dots,K\}$ be the set of companies, $\scr{N}=\{1,2,\dots,N\}$ be the set of consumers, and $\scr{T}=\{1,2,\dots,T\}$ be the finite set of time slots \footnote{\color{blue}{In this paper, we interchangeably use ``time slot" and ``period" to refer to a subdivision of the time horizon.}}.
We formulate a static Stackelberg game between utility companies (the leaders) and their consumers (the followers) to find revenue maximizing prices and optimal demands. In Stackelberg games, the leader(s) first announce their decisions to the follower(s), and then the followers respond. In our game, the leaders send price signals to the consumers, who respond optimally by choosing their demands. To capture the  market competition between the utility companies, we let them play a price-selection Nash game. The equilibrium point of the price-selection game is what utility companies announce to their consumers. The consumers, on the other hand, do not face a game among themselves as they are individual utility maximizers. Figure 1 illustrates the hierarchical interaction between companies and consumers. 
In the parlance of dynamic game theory \cite{basar}, we are dealing here with open-loop information structures, with the corresponding equilibrium at the companies level being open-loop Nash equilibrium. Therefore, this is a one-shot game at which all the prices for all the periods are announced at the beginning of the game, and the followers respond to these prices by solving their local optimization problems. 
\subsection{Consumer-Side} 
Because of energy scheduling and storage, consumers may have some flexibility on when to receive a certain amount of energy. We are concerned about the total amount of shiftable energy. {\color{blue}Period-specific constraints can be added to include non-shiftable energy demand in the problem formulation, as discussed later in Section \ref{generalizations}}. Each energy consumer $n\in \scr{N}$ receives all price signals from each company $k\in \scr{K}$ at each time slot $t\in \scr{T}$ and aims to select his corresponding utility-maximizing demand $d_{n,k}(t)\geq0$ for each time slot from each company, subject to budget and energy need constraints. Denote the price of company $k$ at time $t$ by $p_k(t)$. Let $B_n \geq 0$  and $E^{{\rm min}}_n \geq 0$ denote, respectively, the budget of consumer $n$ and minimum energy need for the entire time-horizon. The utility of consumer $n$ is defined as 
\begin{equation}{ u_n(\mathbf{d}_{n})=\gamma_n\sum_{k\in \scr{K}}\sum_{t\in \scr{T}}\ln(\zeta_n+d_{n,k}(t))} \label{consumer}\end{equation} 
where $\gamma_n>0$ and $\zeta_n\geq1$ are preference parameters.  Note that if $0\leq\zeta_n<1$ or $\gamma_n<0$, the utility of the consumer becomes negative, which is not realistic for demand response applications, and hence we take $\gamma_n>0$ and $\zeta_n\geq1$. A typical value for $\zeta_n$ is $1$, but we still solve the problem for arbitrary $\zeta_n\geq1$ to keep it general. {\color{black}The logarithmic function (\ref{consumer}) is known to provide proportional fairness and is widely used to model consumer behavior in economics \cite{srikant,shadow,gao,basarDR,basarDR2}, {and it has been validated for demand response applications \cite{gao,fan2,han,DRadaptation,sabita}}. Our analysis in this paper is quite general and can be used in any market arrangement with multiple sellers and buyers under budget limitations and capacity constraints.}
Consumer $n$ aims to achieve the highest payoff while meeting the threshold of minimum amount of energy and not exceeding a certain budget. To be more precise,
given $B_n \geq 0$, $E^{{\rm min}}_n \geq 0$, and $p_k(t)>0$, the consumer-side optimization problem is formulated as follows:
\begin{eqnarray}
\underset{\mathbf{d}_{n}}{\hbox{maximize}} && u_n(\mathbf{d}_{n}) \nonumber \\
\hbox{subject to} && \sum_{k\in \scr{K}}\sum_{t\in \scr{T}}p_k(t)d_{n,k}(t)\leq B_n \label{cc1}\\
&& \sum_{k\in \scr{K}}\sum_{t\in \scr{T}}d_{n,k}(t)\geq \,E^{{\rm min}}_n \label{cc}\\
&&d_{n,k}(t)\geq 0,\;\;   \forall k \in \scr{K},\;\;  \forall t \in \scr{T} \label{cc3}
\end{eqnarray}
Note that, as indicated earlier, there is no game played among the consumers. Each consumer responds to the price signals using only her local information.We indirectly handle consumers' cost minimization via our analysis in later sections.
\subsection{Company-Side}
Let the prices chosen by other companies be ${\bf{p_{-k}}}$. The revenue for company $k$ is then given by
\begin{equation}{\pi_k}({\bf{p}}_{k},{\bf{p_{-k}}}):=\sum_{t\in \scr{T}}p_k(t)\sum_{n\in \scr{N}}d_{n,k}({\bf{p}}_k,{\bf{p_{-k}}},t).\label{UC}\end{equation}
Given the power availability of company $k$ at period $t$, denoted by $G_k(t)$, and for a fixed ${\bf{p_{-k}}}$, company $k$ solves the following problem:
\begin{eqnarray}
\underset{\mathbf{p}_k}{\hbox{maximize}} && \pi_k ({\bf{p}}_k,{\bf{p}_{-k}})
\nonumber \\
\hbox{subject to} && \sum_{n\in \scr{N}} d_{n,k}({\bf{p}}_k,{\bf{p}_{-k}},t) \leq G_k(t),\;\; \forall \; t \in \scr{T} \label{prob2} \\
&& p_k(t)> 0, \;\; \forall \; t \in \scr{T} \label{end}
\end{eqnarray}
  The goal of each company is to maximize its revenue\footnote{In later sections we show how companies can alter their problems to profit-maximization instead}. Additionally, because of the market competition, the prices announced by other companies also affect the determination of the price at company $k$. Thus, company $k$ selects its price in response to what other competitors in the market have announced; this response is also constrained by the availability of power. Thus, what we have is a Nash game among the companies. We emphasize that while each company's problem is affected by what its competitors decide, we can still achieve the equilibrium strategies using only local information, via our distributed algorithm discussed later in Section \ref{algorithm}. Finally, while at this point we have ${\bf G}_k$ fixed for each company $k$, we will later formulate a power allocation game to optimally choose them.

\section{Demand Selection and Revenue Maximization (Stackelberg Game)} \label{game1}
In this section, we solve the above optimization problems in closed form and show that the solutions are unique.


\subsection{Consumer-Side Analysis}
We start by relaxing the minimum energy constraint (\ref{cc}).
  For each consumer $n\in \scr{N}$, the associated Lagrange function is given as follows:
\begin{eqnarray*}
L_n &=& \gamma_n\sum_{k\in \scr{K}}\sum_{t\in \scr{T}}\ln(\zeta_n+d_{n,k}(t))
 \\
&&- \lambda_{n,1}\left(\sum_{k\in \scr{K}}\sum_{t\in \scr{T}}p_k(t)d_{n,k}(t)-B_n\right) \\
&&+\sum_{k\in\scr{K}} \sum_{t\in\scr{T}} \lambda_{n,2}(k,t)d_{n,k}(t)
\end{eqnarray*}
  where ${\bf \lambda_{n}}$ are the Lagrange multipliers. The KKT conditions of optimality in this case are sufficient because the objective function is strictly concave and the constraints are linear \cite{NL}, and solving for them leads to \begin{equation}
d^*_{n,k}(t)= \frac{B_n+\sum_{j\in \scr{K}}\sum_{h\in \scr{T}}p_j(h)\zeta_n}{KTp_k(t)}-\zeta_n,  \; \forall \;t\in \scr{T}, \; k\in \scr{K},  \label{xx}
\end{equation}
which is a generalization of the single-period case in \cite{sabita}. A detailed derivation of (\ref{xx}) can be found in \cite{mywork}. {\color{black} We remark that $d^*_{n,k}(t) \geq 0$  {\color{blue} because the objective function is strictly increasing.}}

The following theorem, whose proof can be found in the Appendix, states the necessary and sufficient condition for $B_n$ so that the above demands meet the minimum energy constraint (\ref{cc}).

\begin{theorem}
For each consumer $n \in \scr{N}$, the demands $d^*_{n,k}(t)$ given by (\ref{xx}) satisfy (\ref{cc}) if, and only if, 
\begin{equation} B_n \geq \frac{E_n^{{\rm min}}+\zeta_nKT}{\sum_{k\in \scr{K}}\sum_{t\in \scr{T}}\frac{1}{KTp_k(t)}}-\zeta_n \sum_{k\in \scr{K}}\sum_{t\in \scr{T}}p_k(t). \label{budget} \end{equation}
\end{theorem}

{\color{black}
\begin{remark} The above theorem can be interpreted as billing costs minimization. At the equality of (\ref{budget}), $B_n$ corresponds to the minimum budget needed for consumer $n$ to satisfy his energy need constraint, given the set of prices chosen by utility companies. Such a minimum $B_n$ can serve as a theoretical benchmark in which one can measure whether or not consumers are paying more than what is necessary.  We later demonstrate that with real data from demand response experiments, using the equality in (\ref{budget}) leads to savings in the range of $10\%-30\%$. \hfill $\Box$
\end{remark}

\begin{assumption}
For each consumer $n$, the budget $B_n$ satisfies the condition (\ref{budget}).
\label{assumption1}
\end{assumption}
}
\subsection{Company-Side Analysis} 
We apply the demands derived in the consumers-side analysis (which were functions of the prices) and show that optimality is achieved at the equality of the constraint (\ref{prob2}). We start by solving for prices that satisfy the equality at (\ref{prob2}) and then prove that they are revenue-maximizing, strictly positive, and unique. 
Consider the equality in (\ref{prob2}), and by the optimal demands (\ref{xx}), there holds
$$
\frac{\sum_{n\in \scr{N}}B_n+\sum_{n\in \scr{N}}\zeta_n\sum_{j\in \scr{K}}\sum_{h\in \scr{T}}p_j(h)}{KTp_k(t)} =\sum_{n\in \scr{N}}\zeta_n + G_k(t),
$$
for all $t \in \scr{T}$. 
Let $Z=\sum_{n\in \scr{N}}\zeta_n$ and $B=\sum_{n\in \scr{N}}B_n$. Then, for each company $k \in \scr{K}$,
\begin{equation} B+Z\sum_{j\in \scr{K}}\sum_{h\in \scr{T}}p_j(h) = KTp_k(t)(G_k(t)+Z),\;\; \forall \;t \in \scr{T}.  \label{AP} \end{equation}
The above equation (\ref{AP})  can be presented as the following system of linear equations
\begin{equation}AP=Y,\label{AP2}\end{equation}
%
where $A$ is a $KT\times KT$ matrix 
whose diagonal entries are $KT(G_k(t)+Z)-Z$, $k\in\scr{K}$, $t\in\scr{T}$,
and off-diagonal entries all equal to $-Z$, 
$P$ is a vector in $\R^{KT}$ stacking $p_k(t)$, $k\in\scr{K}$, $t\in\scr{T}$,
and $Y$ a vector in $\R^{KT}$ whose entries all equal to $B$.

We have the following results (proofs are in the Appendix).
\begin{lemma}
The matrix $A$ is invertible.
\end{lemma}

\begin{lemma}
The prices that solve (\ref{AP2}) are strictly positive and are unique. For each $t\in\mathcal{T}$, $k\in\mathcal{K}$, the price is given by
    \begin{equation} p^*_k(t)=\frac{B}{G_k(t)+Z}\left(\frac{1}{KT-\sum_{j\in \scr{K}}\sum_{h\in \scr{T}}\frac{Z}{G_j(h)+Z}}\right),\label{p}\end{equation}
where $B=\sum_{n\in \scr{N}}B_n$ and $Z=\sum_{n\in \scr{N}}\zeta_n$.
\end{lemma}

\begin{remark} 
Letting $\zeta_n=1$ for each consumer, the value of $Z$ coincides with $N$. In this case, by (\ref{p}), we observe that for any given ${\bf G_k}$, the price $p^*_k(t)(G_k(t)+N)$ is a constant for all $t \in \scr{T}$ and $k \in \scr{K}$.
Thus, the power availability is inversely proportional to the prices. \hfill$\Box$
\end{remark}

\begin{remark} 
Lemma 2 provides a computationally cheap expression for the prices. Since $p^*_k(t)$ can be directly computed using (\ref{p}), there is no need to numerically compute $A^{-1}$ or $|A|$ to solve (\ref{AP2}). This enables us to deal with a large number of periods or utility companies, without worrying about computational complexity.\hfill$\Box$
\end{remark}

{\color{black}Due to production costs and market regulations, $p^*_k(t)$ cannot be outside the range of some lower and upper bounds $[p^{{\rm min}}_k(t),p^{{\rm max}}_k(t)]$  for all $t \in \scr{T}$ and $k \in \scr{K}$, as in \cite{sabita}. If $p^*_k(t)<p^{{\rm min}}_k(t)$, then $p^*_k(t)$ is set to $p^{{\rm min}}_k(t)$, and similarly for the upper-bound, if  $p^*_k(t)>p^{{\rm max}}_k(t)$, then we set  $p^*_k(t)=p^{{\rm max}}_k(t)$. Accordingly, denote the strategy space of utility company $k$ (a leader in the game) at $t$ by $\scr{L}_{k,t}:=[p^{{\rm min}}_k(t),p^{{\rm max}}_k(t)]$. The strategy space of $k$ for the entire time horizon is $\scr{L}_{k}=\scr{L}_{k,1}\times\dots\times\scr{L}_{k,T}$.
 The strategy space of all companies is $\scr{L}=\scr{L}_{1}\times\dots\times\scr{L}_{K}$. {\color{blue} For given price selections  ${\bf{p}}:=({\bf{p}}_1,\dots,{\bf{p}}_K) \in \scr{L}$}, the optimal response from all consumers is
$${\bf{d^*(p)}}=\{{\bf{d}}_1^*({\bf{p}}),{\bf{d}}_2^*({\bf{p}}),\dots,{\bf{d}}_N^*({\bf{p}})\}$$
where for each $n \in \scr{N}$, ${\bf{d}}^{*}_{n}({\bf{p}})$ is the unique maximizer for $u_n({\bf{d}}_n,{\bf{p}})$ and is given by (\ref{xx}).
}
We now have the following theorem, whose proof can be found in the Appendix.

\begin{theorem}[Existence and Uniqueness of the Stackelberg Equilibrium]Under Assumption \ref{assumption1}, the following statements hold: 
\begin{itemize}
\item[(i)] There exists a unique (open-loop) Nash equilibrium for the price-selection game and it is given by (\ref{p}).\\
\item[(ii)] There exists a unique (open-loop) Stackelberg equilibrium, and it is given by the demands in (\ref{xx}) and the prices in (\ref{p}).
\end{itemize}
\label{mainTHM}
\end{theorem}

At the Stackelberg equilibrium, it can easily be verified that  
\begin{equation}\sum_{k\in \scr{K}}\pi_k({\bf{p}}^*_k,{\bf{p}^*_{-k}})=\sum_{n\in \scr{N}}B_n.\label{eq1}\end{equation}
One observation is that when a company gains in terms of revenue, the same amount must be lost by other companies because the sum of revenues is a constant, which demonstrates a conflict of objectives between utility companies. However, by the definition of the equilibrium strategy, this is the best each company can do, for fixed power availabilities ${\bf G_k}$. But, given a total amount of available power, $G^{{\rm total}}_k$, a company has across the time horizon, it is possible that it gains in terms of revenue by an efficient power allocation. This motivates us to formulate a power allocation game and analytically answer the following question: How can company $k$ allocate its power so that it maximizes its revenue? {\color{black} Furthermore, for now, for ease of exposition, we neglect network and other company-specific constraints. Such considerations are later discussed in Section \ref{generalizations}.  For the remaining part of this paper, unless otherwise stated, we also have the following simplifying assumption.

\begin{assumption}
For each consumer $n$, we have $$\gamma_n=\zeta_n=1.$$ \label{assumption2}
\end{assumption}

The above assumption implies that $Z$ is equal to the number of consumers $N$.}

\section{Power Allocation (Nash Game)} \label{game2}

\begin{figure*}
\centering
\includegraphics[width=.9\linewidth]{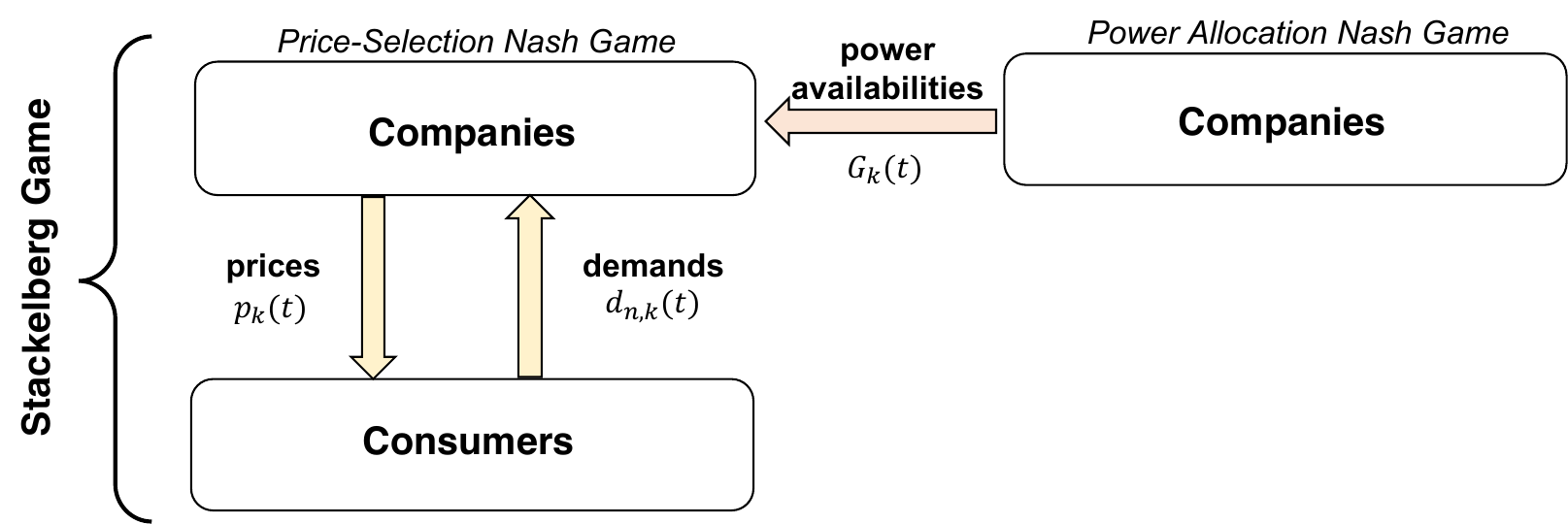}
\caption{The interaction between companies and their consumers, along with power allocation. First, companies play a Nash power allocation game. Once power availabilities are allocated across all periods, companies and consumers play the Stackelberg game which dictates optimal prices and demand selection.}
\label{sketch_DR2}
\end{figure*}
In this section, we exploit the closed-form solutions for consumer demands and companies' prices to formulate and solve a power allocation game for companies. We note that while we use the closed-form solutions to define the power allocation game, it is to be played {\em{before}} the Stackelberg game, and its outcomes define the fixed power availabilities in the constraints of the companies in the Stackelberg game. Given the power availabilities from other companies, ${\bf{G_{-k}}}$, and since the equality in (\ref{prob2}) is satisfied at equilibrium, the revenue function of company $k$ can be represented as
\begin{equation} \pi_k({\bf{G}}_k,{\bf{G_{-k}}})=\sum_{t\in \scr{T}}p^*_k(t)G_k(t). \end{equation}  The optimal prices (\ref{p}) are functions of ${\bf{G}}_k$ and ${\bf{G_{-k}}}$, leading to the revenue function being equal to
\begin{equation} B\sum_{t\in \scr{T}}\frac{G_k(t)}{(G_k(t)+N)(KT-\sum_{j\in \scr{K}}\sum_{h\in \scr{T}}\frac{N}{G_j(h)+N})}, \label{Uk} \end{equation}
where $B=\sum_{n\in \scr{N}}B_n$. {Note that company $k$ receives a {\it fraction} of the total budgets. This fraction depends on what company $k$ offers in the multi-period-multi-company demand response framework, and what other companies also offer. Thus, when company $k$ can change what it offers, it can potentially increase the fraction it receives, and the power allocation game becomes natural, since the revenue function depends on other players' decisions. } For this game, which can be played before the Stackelberg game, which we have already solved, companies allocate their powers across all periods, and the outcome dictates the fixed power availabilities for the Stackelberg game. Figure \ref{sketch_DR2} provides an illustration. 



Let the total capacity for company $k$ for the entire time horizon be $G^{{\rm total}}_k$. Denote the action set of company $k$ at time $t$ by $\scr{P}_{k,t}:=[0,G^{{\rm total}}_k]$. Thus, given ${\bf{G_{-k}}}$,  the company $k$ solves the following problem:
\begin{eqnarray}
\underset{\mathbf{G}_k }{\hbox{maximize}} && \pi_k({\bf{G}}_k,{\bf{G_{-k}}})
\nonumber \\
\hbox{subject to} && \sum_{t\in \scr{T}}G_k(t) \leq G^{{\rm total}}_k, \\ \label{prob3}
&& G_k(t)\geq0, \ \forall t\in\scr{T}. \nonumber\end{eqnarray}

{\color{black}The above problem is only applicable for the case when generation is fully controllable. For the smart grid, because of the availability of various generation sources, full-controllability does not always hold, and in fact, for renewable resources it could be completely gone. We demonstrate the possibility of relaxing this assumption later in Section \ref{generalizations}.}

\subsection{Existence and Uniqueness of Nash Equilibrium} 
%
The following theorem, whose proof can be found in the Appendix, states the existence and uniqueness of Nash equilibrium in the power allocation game, and provides an expression for it.
\begin{theorem}
Under Assumptions \ref{assumption1}-\ref{assumption2}, if ${\bf G}_k$ is fully controllable, there exists a unique pure-strategy Nash equilibrium for the power allocation game,
and it is given by \begin{equation} G^*_k(t)=\frac{G^{{\rm total}}_k}{T} \,\,\,,\,\, \forall \,t\in \scr{T}, \forall \, k\in \scr{K}.\label{NE}\end{equation}\label{allocation}\end{theorem}

Interestingly, the optimal strategy for each company is to equally allocate its power across all time periods. The proof of Theorem \ref{allocation} reveals that (\ref{Uk}) is strictly concave and increasing in each $G_k(t)$. This is an important property that allows accommodating further company-specific operational constraints and relaxing the full-controllability assumption. To illustrate, suppose that company $k$ has a mix of generation sources for which generation is controllable for some periods and only partially controllable for others. Then, it can add linear constraints to problem (\ref{prob3}) reflecting inter-temporal considerations at the generation-side (such as ramping limits). Existence and uniqueness of a pure-strategy Nash equilibrium are still guaranteed due to the strict concavity of the objective \cite{basar}. Since generation costs are typically assumed to be convex \cite{bosebook} (denote it by $c_k$ for each company $k$), company $k$ can also allocate its generation to maximize its profit, by subtracting the cost from (\ref{Uk}). One can alter the objective function of the power allocation game to 
   \begin{align}&B\sum_{t\in \scr{T}}\frac{G_k(t)}{(G_k(t)+N)(KT-\sum_{j\in \scr{K}}\sum_{h\in \scr{T}}\frac{N}{G_j(h)+N})}\nonumber \\
   & \qquad\qquad\qquad \qquad-\sum_{t\in \scr{T}} c_k(G_k(t)),\label{cvx}\end{align}
   
 \noindent and the problem reflects profit-maximization in this case. Using (\ref{cvx}) and following our analysis, {\color{blue} we conclude that each company maximizes a strictly concave function}, and one can easily conclude the existence of a pure-strategy Nash equilibrium in this case. 

\section{{\color{black}Distributed Algorithm} }\label{algorithm}
The Nash equilibrium (NE) for the power allocation game given by (\ref{NE}) can easily be computed by each company $k$ using its local information. Moreover, for consumers, it can be seen from (\ref{xx}) that in the computation of optimal demand selection for consumer $n$, no information from other consumers is needed, and consumer $n$ only needs local information for optimal response. However, the closed-form solution for optimal prices given by (\ref{p}) requires each company $k$ to know consumers' budgets and the power availability of all the other companies. Companies might not want to share such information with each other. To circumvent such a privacy concern, we propose a distributed algorithm that allows companies to compute their optimal prices using only local information, and show that this algorithm converges to the optimal prices given by  (\ref{p}). The algorithm, combined with utility-maximizing demands given by (\ref{xx}) and the NE given by (\ref{NE}), leads to the computation of all the optimal strategies with only local information at both the company level and the consumer level. 
\begin{algorithm}
\begin{algorithmic}[1]
\State Arbitrarily choose $p^{(0)}_k(t) ,\,\,\forall t\in\scr{T}, \,\,\,\forall k\in\scr{K}$ 
\State Repeat for $i=1, 2,3,\dots$ \label{next}
\State For each consumer $n \in \scr{N}$, compute $d^{(i)}_{n,k}(t)$ from  $k\in\scr{K}$ at  $t\in\scr{T}$ by (\ref{xx}), then update utility companies with demand signals\label{user}
\State {\color{black} Pick a company $k\in\mathcal{K}$ at time $t\in\mathcal{T}$ such that $p^{(i+1)}_k(t)$ is not yet computed, and compute it using (\ref{update}) \label{UCupdate}}
\State If $p^{(i+1)}_k(t)\neq p^{(i)}_k(t)$, update consumers and go to \ref{user}
\State Else, send a no-change signal to consumers and go to \ref{UCupdate}
\State If $p^{(i+1)}_k(t) = p^{(i)}_k(t)\,
\,\forall t\in\scr{T}, \,\,\,\forall k\in\scr{K}$, stop 
\State Else, go to \ref{next}
\end{algorithmic}
\caption{Distributed algorithm for computing the prices with local information} 
\end{algorithm}

For each iteration $i\in\{0,1,2,\ldots\}$, denote the demand from consumer $n$ at time $t$ from company $k$ by $d^{(i)}_{n,k}(t)$, and the price announced by company $k$ and time $t$ by $p^{(i)}_k(t)$. In our algorithm,  $p^{(0)}_k(t)$ is chosen arbitrarily for each company $k\in \scr{K}$ and time $t\in\scr{T}$. Based on the initial price selection, $d^{(0)}_{n,k}$ is computed using (\ref{xx}).  Then, the prices are sequentially updated using the following update rule: 
\begin{equation}p^{(i+1)}_k(t)=p^{(i)}_k(t)+\frac{\sum_{n\in \scr{N}}d^{(i)}_{n,k}(t)-G_k(t)}{\epsilon^{(i)}_{k,t}},\label{update}\end{equation}
where $\epsilon^{(i)}_{k,t}>0$ is appropriately selected for company $k$ at time $t$ in iteration $i$, and we present an expression for it as a function of $p^{(i)}_k(t)$ in Theorem \ref{d_algorithm}. 
Whenever a company $k$ updates its price at time $t$, it transmits the price to each consumer $n\in\scr{N}$, and they modify their demands accordingly. Once prices converge to their optimal values, consumers optimally respond by (\ref{xx}) and the algorithm terminates. We have the following theorem for the convergence of the algorithm; its proof can be found in the Appendix.

 \begin{theorem}
Under Assumptions \ref{assumption1}-\ref{assumption2}, for each company $k\in\scr{K}$ at time $t\in\scr{T}$ in iteration $i\in\{0,1,2,\ldots\}$, if the prices are sequentially updated using (\ref{update}) such that
$$ \epsilon^{(i)}_{k,t} =  \frac{G_k(t)+N}{p^{(i)}_k(t)} + \delta,$$
where $\delta\geq0$, then Algorithm 1 converges to optimal prices.
\label{d_algorithm}
 \end{theorem}

\section{Asymptotic Regimes} \label{asymp}
{\color{black}In this section, we study the asymptotic (limiting) behavior as $T\rightarrow \infty$ or $N\rightarrow \infty$. While neither $T$ or $N$ can be arbitrarily large in practice, analyzing the asymptotic behavior brings in deep insights. For example, it reveals that consumers benefit as $T$ grows. As $N$ grows, our asymptotic analysis allows us to compute an appropriate company-to-consumer ratio $\frac{K}{N}$. We show these insights by studying how the utility functions, revenues, prices, and demands are affected as $T$ or $N$ grows. {For the rest of this section, in addition to Assumptions \ref{assumption1}-\ref{assumption2}, we assume the following. 
\begin{assumption} 
The total power available for the entire time horizon $G_k^{{\rm total}}$ is the same for each company $k\in\scr{K}$.
\label{assumption3}
\end{assumption}
}}
\subsection{When the Number of Periods Grows} Under Assumptions \ref{assumption1}-\ref{assumption3}, at equilibrium, it follows that 
optimal prices and demands are given by 
\begin{equation} p^*_k(t)=\frac{\sum_{m\in \scr{N}}B_m}{KTG^*_k(t)}=\frac{\sum_{m\in \scr{N}}B_m}{KG^{{\rm total}}_k}, \label{ppp} \end{equation}
\begin{equation}d^*_{n,k}(t)=\frac{B_n+KTp^*_k(t)}{KTp^*_k(t)}-1=\frac{G^{{\rm total}}_kB_n}{T\sum_{m\in \scr{N}}B_m},\,\, \label{t1} \end{equation}
and the utility of consumer $n$ becomes\begin{equation} u_n=KT\ln\left(1+\frac{G^{{\rm total}}_kB_n/\sum_{m\in \scr{N}}B_m}{T}\right),\end{equation}
in which $G^{{\rm total}}_kB_n/\sum_{m\in \scr{N}}B_m$ is positive. Thus, as $T$ increases, the multiplicative term $ KT$ of the logarithmic function increases at a faster rate than the decrease of 
$\ln\left(1+{B_nG^{{\rm total}}_k/B}/{T}\right)$. 
Hence, as $T$ increases, the utility of each consumer $n \in \scr{N}$ monotonically increases.
Taking the limit, it can be verified that 
\begin{equation}\lim_{T\rightarrow\infty} u_n(T)=\frac{KG^{{\rm total}}_kB_n}{\sum_{m\in \scr{N}}B_m}.\end{equation}
Furthermore, note that the demand $d^*_{n,k}(t)$ from consumer $n\in \scr{N}$ from company $k\in \scr{K}$ at time $t\in \scr{T}$ converges to zero as $T \rightarrow \infty$. We claim that the revenues are constants. To see this, recall that
\begin{align*}
\pi_k({\bf{p}}^*_k,{\bf{p}^*_{-k}}) &= p^*_k(t)G^{{\rm total}}_k = \frac{\sum_{m\in \scr{N}}B_m}{K},
\end{align*}
which is a constant since both the number of companies and the budgets of the consumers are fixed.

\begin{remark} At the equilibrium, the monotonicity of the utilities of the consumers shows that increasing the number of periods leads to more incentives for consumers' participation in demand response. However, it might not be very beneficial to increase the number of periods to a very high value. First, the rate of increase in terms of consumers' utilities gets progressively smaller. Second, having a high number of periods leads to smaller demands for each period and that might violate some minimum energy need for particular periods at the consumers' level. So, it is beneficial to increase the number of periods up to a certain point (compared to having $T=1$), but it might not be beneficial to let $T$ become arbitrarily large. 
\hfill$\Box$
\end{remark}

\begin{remark} Note that the limit point of the utility function of consumer $n$ is the proportion of his budget to the total budgets times the total power availability. So if a particular consumer has $1\%$ of the sum of all the budgets, he gets $1\%$ of the available power. Furthermore, the revenue for each company is the proportion of the sum of the budgets to the number of companies. In addition, the demand by consumer $n$ from company $k$ at time $t$ is the proportion of his budget to the total budgets times the total power availability at $t$ from $k$.
\hfill$\Box$
\end{remark}

\subsection{When the Number of Consumers Grows} When the number of consumers increases, each additional consumer has some budget $B_n$. With the total power availability from companies being fixed, they will increase their prices. {\color{black}  We have the following simplifying assumption. 

\begin{assumption} 
The budget for each consumer $n \in \scr{N}$ is the same.
\label{assumption4}
\end{assumption}

Under Assumptions \ref{assumption1}-\ref{assumption4}, we increase the number of consumers $N$ and see what happens as $N \rightarrow \infty$.} In this case, the optimal prices and demands become
\begin{eqnarray}p^*_k(t)&=&\frac{NB_n}{KTG^*_k(t)} \label{pppp}\\
d^*_{n,k}(t) &=& \frac{G^{{\rm total}}_k}{TN} \label{dd} \end{eqnarray}
Clearly, $p^*_k(t)\rightarrow\infty$ as $N\rightarrow\infty$ and $d^*_{n,k}(t)\rightarrow0$ as $N\rightarrow\infty$. When the population is large and the power availability is fixed, it is not surprising that $d^*_{n,k}(t)\rightarrow0$ because the portion each consumer can get from the available power gets smaller and smaller as $N$ increases. Furthermore, it can be easily verified that $\lim_{N\rightarrow \infty}\pi_k(N)=\infty$ and $\lim_{N\rightarrow \infty}u_n(N)=0$. Thus, with the limit points resulting in unrealistic outcomes, a balance between the supply and demand needs to be achieved, which we do by finding an appropriate company-to-consumer ratio. 

Now, the question we ask is: For a given maximum allowable market price $p^{{\rm max}}_k(t)$, call it $p^{{\rm max}}$, what is the appropriate company-to-consumer ratio $\frac{K}{N}$? If there are more companies than necessary in the market, there will be losses in terms of revenues. On the other hand, if there are fewer companies than necessary, the prices can exceed $p^{{\rm max}}$, leading to undesirable outcomes. The following theorem, whose proof can be found in the Appendix, provides an optimal ratio at which prices do not exceed $p^{{\rm max}}$ and the revenues being maximized while satisfying the equality in (\ref{eq1}). 
{\color{black}
\begin{theorem}
Under Assumptions \ref{assumption1}-\ref{assumption4}, at the NE of the power allocation game, and at the Stackelberg equilibrium of the price and demand selection game, the optimal prices given by (\ref{p}) satisfy  
\begin{eqnarray*}
 p^*_k(t) &\leq &p^{{\rm max}},\\
 \sum_{k\in \scr{K}}\pi_k({\bf{p}}^*_k,{\bf{p^*_{-k}}})&=&\sum_{n\in \scr{N}}B_n,
\end{eqnarray*}
if, and only if, \begin{equation*} \frac{K}{N} \geq \frac{B_n}{p^{{\rm max}}TG^*_k(t)}, \end{equation*}
for each $t \in \mathcal{T}$ and $k\in\mathcal{K}$.
\end{theorem}}


\section{Case Studies} \label{numerical}

In this section, we present results on some case studies on representative days from a Dutch smart grid pilot \cite{dutch} and the EcoGrid EU project \cite{ecogrid}. We numerically study optimal prices and demands, and their corresponding payments and utility functions. We show how our approach results in monetary savings for consumers. Furthermore, we show that increasing the number of periods provides more incentives for consumers' participation in demand response management. Additionally, we demonstrate the fast convergence of our distributed algorithm to optimal prices. We also release an open-source interactive tool containing the simulations in \cite{tool2}. For both the Dutch smart pilot and the EcoGrid EU projects, the data are unavailable in raw format. Thus, whenever it is needed, we estimate some data points from figures available in the corresponding references \cite{dutch,ecogrid}. 

Recall that at the Stackelberg equilibrium, the total power availabilities ${\bf G}$ match the aggregate demands. That is, 
$$ \sum_{n \in \mathcal{N}} d^*_{n,k}(t)=G_k(t), \qquad \forall t\in\mathcal{T}, k\in\mathcal{K}.$$ Here, we use the experimental hourly variation of the total demands to choose values for ${\bf G}$ and the minimum energy need ${\bf E}^{\min}$. This allows us to establish a common aspect between our results and the experimental results, so that we can appropriately explore how our framework compares to real-life experiments. We also use the lower-bound on the minimum budget condition (\ref{budget}), so that we can also quantify potential savings.  
{\color{blue}  From the consumers' perspective, the prices are given parameters in both our model and the experimental setups. The optimal demands are functions of the prices, and the optimal prices naturally depend on the parameters of the consumers and companies. To bring deep insights, we make the differentiating aspect between our model and the experimental results an economic one. And hence, we pick the parameters such that the equilibrium demands and experimental ones are similar, but the prices, and essentially what consumers pay, are different. Utilizing Theorem 1, we conclude that the equilibrium prices bring savings to consumers, and by definition, they automatically consider the incentives of companies as they are revenue-maximizing. A main conclusion of this paper is that this quantifies the economic gap, in terms of consumer savings, between our game-theoretic benchmark and existing experimental results. On the other hand, in our analysis, we have relaxed some constraints for tractability, such as power flow and demand inelasticity considerations, and it remains open to explore the underlying tradeoffs, since adding these considerations might reduce the potential monetary savings for consumers. Such considerations were not directly included in the models studied in \cite{ecogrid,dutch}, as their focus was to experimentally explore the consumers' behavior in response to changing prices. It is worth mentioning that the results in \cite{dutch} revealed that consumers are mainly flexible about adjusting the consumption of white goods (washing machine, dishwasher, etc). It was also concluded in \cite{ecogrid} that demand response did not result in distribution feeder congestion relief, and consumers with automatic equipment were the most responsive ones. Nevertheless, later in Section IX, we demonstrate how our framework can be utilized to include additional network and consumer-specific and/or company-specific constraints, which make it possible to add constraints for congestion relief. Including such constraints will likely make it necessary to compute the equilibrium prices and demands algorithmically, which we relegate to future endeavors, as we emphasize here more on revealing deep insights via having tractable analysis.}
 \begin{figure*}
\centering
\includegraphics[width=\linewidth, height=2in]{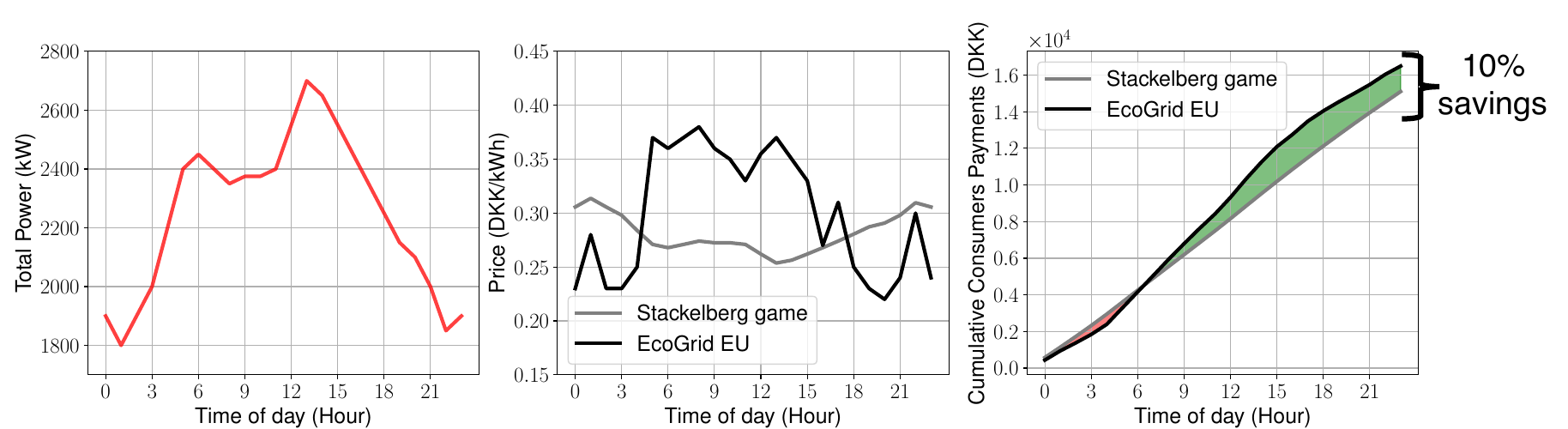}
\caption{Total power offered by company (left), Stackelberg game and EcoGrid EU experimental prices (middle), and the cumulative payments and billing savings for all consumers (right). }
\label{EcoGridEU_Price_Power}
\end{figure*}

\subsection{EcoGrid EU Project} This demand response project was conducted from March 2011 to August 2015 in Bornholm, Denmark. The number of consumers in this experiment was approximately $2000$. For a representative day (December 5th, 2014), we apply our method to hourly prices and shiftable demand consumption from this experiment. The experimental prices are in  \text{DKK}/MWh and we scale them to  \text{DKK}/kWh. We start by assuming that there is only one company ($K=1$) and letting the consumers to be homogeneous (they have the same budgets and energy need) with $N=2000$, and then generalize the results to $K>1$ and heterogeneous consumers. Since we are taking hourly prices for a day, we have $T=24$. 

\subsubsection{Finding the necessary parameters}In our model, for each period $t$, we have a fixed power availability $G_1(t)$ on the supply-side. Also, for each consumer $n$, his minimum demand $E^{\min}_n$ and budget $B_n$ are fixed for the entire horizon. These are necessary parameters that need to be known to solve for optimal demands and prices. We let the power availabilities ${\bf G}_1$ match the experimental hourly variation of the total demand. For the entire time-horizon, we have  $$\sum^{2000}_{n=1} E^{\min}_n=\sum^{24}_{t=1}G_1(t) \approx 54 \ \text{{\color{blue}MWh}}.$$ For homogenous consumers, it follows that $$E^{\min}_n= \frac{\sum^{24}_{t=1}G_1(t)}{2000} \approx 27 \ \text{{\color{blue}kWh}}.$$ Next, using Theorem 1, we plug-in $E^{\min}_n$ and the experimental hourly prices in (\ref{budget}) to find the minimum budget need, which is $B_n\approx 7.6 \ \text{ \text{DKK}}$, for each $n$. 
   
 \subsubsection{Numerical Results} Now, using the parameters found above, we can compute the optimal demands and prices for the Stackelberg game using (\ref{xx}) and (\ref{p}), and study their effects. 
 
In Figure \ref{EcoGridEU_Price_Power}, we plot the total power availabilities ${\bf G}_1$, the prices found experimentally and using the Stackelberg game, and the corresponding total payments by all consumers for their demands. Our approach leads to prices that have a slightly smaller mean than in the experiment and a significantly smaller variance, which is a desirable property \cite{IFAC}. {\color{black} At the equilibrium point, as stated in Remark 2, we observe that $$p^*_k(t)(G_k(t)+N)=p^*_k(t)\Bigg(\sum_{n\in\mathcal{N}}  d^*_{n,k}(t)+N\Bigg)$$  is a constant for each period $t$ and each company $k$. Hence, whenever company $k$ at time $t$ has a large amount of power available to sell $G_k(t)$, it would lower its price, and vice versa. Here, consumers are attracted to buy more whenever the price is low, and will buy less whenever the price is high, which is intuitive.} One advantage our approach has is that it results in billing savings for consumers, as we show in Figure \ref{EcoGridEU_Price_Power} (this demonstrates the importance of Theorem 1, which we use to find the minimum budget need for the consumers). {\color{blue} Here, the equilibrium demands are similar to the experimental values, but since the prices differ,  consumers receive the same amount of energy at smaller costs}. This would lead to more monetary incentives for active consumer participation in demand response management, while being consistent with the company's objectives, since the Stackelberg game prices found using (\ref{p}) are revenue-maximizing as shown in the proof of Theorem 2.
\begin{figure*}
\centering
\includegraphics[width=\linewidth, height=4in]{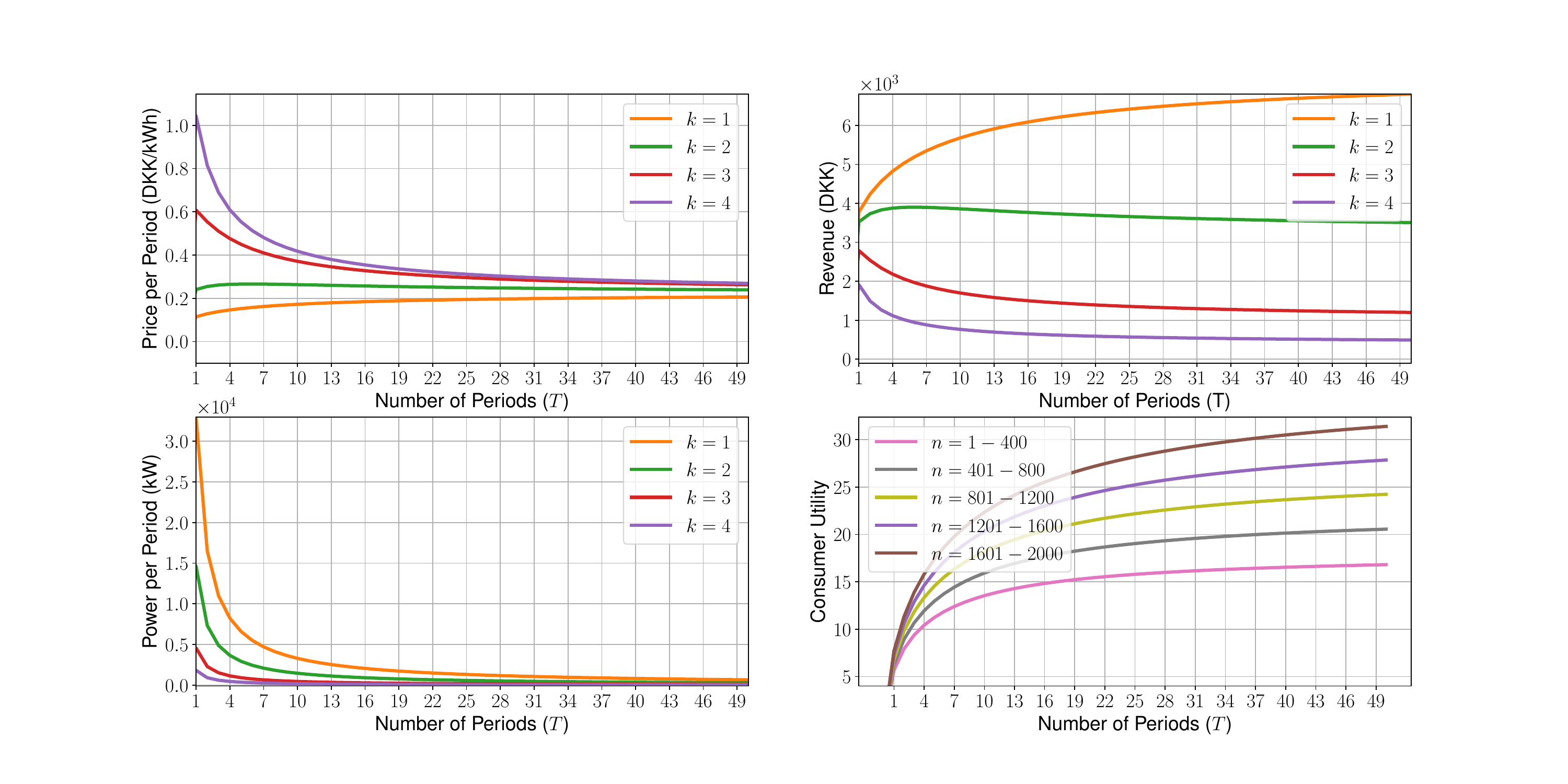}
\caption{The effects of varying the number of periods for companies (with different market shares and at Nash equilibrium of the power allocation game) and heterogeneous consumers (with different budgets) using the EcoGrid EU experimental data.}
\label{Influence_of_T}
\end{figure*}
\begin{figure*}
\centering
\includegraphics[width=\linewidth, height=2in]{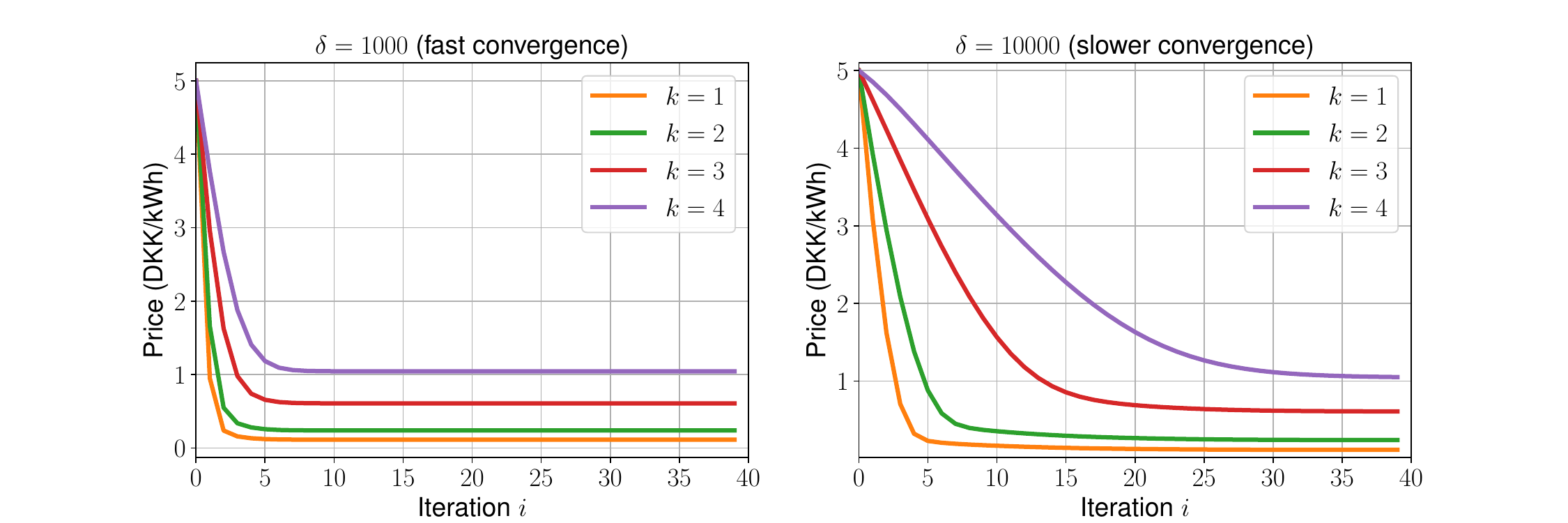}
\caption{ Distributed algorithm's performance (Theorem 4 requires $\delta\geq0$) using the EcoGrid EU experimental data.}
\label{algorithm1}
\end{figure*}

Next, we make consumers heterogeneous and increase the number of companies. We differentiate between consumers by varying their budgets, and take 5 classes of consumers, as in the EcoGrid EU experiment. We let consumers' budgets be $B_{1-400}=4  \ \text{DKK}$, $B_{401-800}=5  \ \text{DKK}$, $B_{801-1200}=6 \  \text{DKK}$, $B_{1201-1600}=7  \ \text{DKK}$, and $B_{1601-2000}=8 \ \text{DKK}$. We also let the number of companies be $K=4$, which is consistent with the actual energy sources used in the experiment. Precisely, the system is powered by 61\% wind energy ($k=1$), 27\% biomass ($k=2$), 9\% solar energy ($k=3$), and 3\% biogas ($k=4$). We split the {\color{black}total need ($54 \ \text{MWh}$)} among the energy sources according to experimental proportions, assuming that each energy source is owned by a single company that acts as a company in our game. 
 
With the above setup, we study the effect of varying the number of periods $T$ from $1$ to $50$. To do this, we need to find a way for companies to allocate their total power across the time horizon for each fixed $T$, which can be done by using Theorem 3, which states that equally splitting the total power across the time horizon for each company $k$ constitutes a unique Nash equilibrium for the power allocation game (it is also shown to be the global maximizer in the proof). 
 \begin{figure*}
\centering
\includegraphics[width=\linewidth, height=2in]{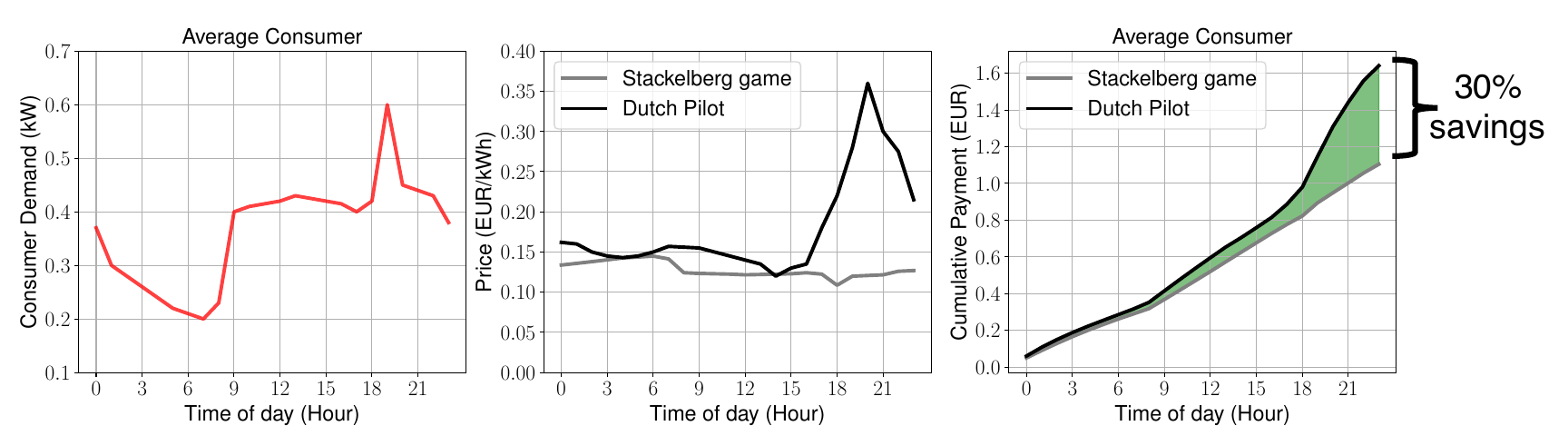}
\caption{Average consumer demand (left), Stackelberg game and Dutch pilot prices (middle), and the cumulative payments and billing savings for average consumer (right).}
\label{Netherlands_Price_Power}
\end{figure*}
\begin{figure*}
\centering
\includegraphics[width=\linewidth, height=2in]{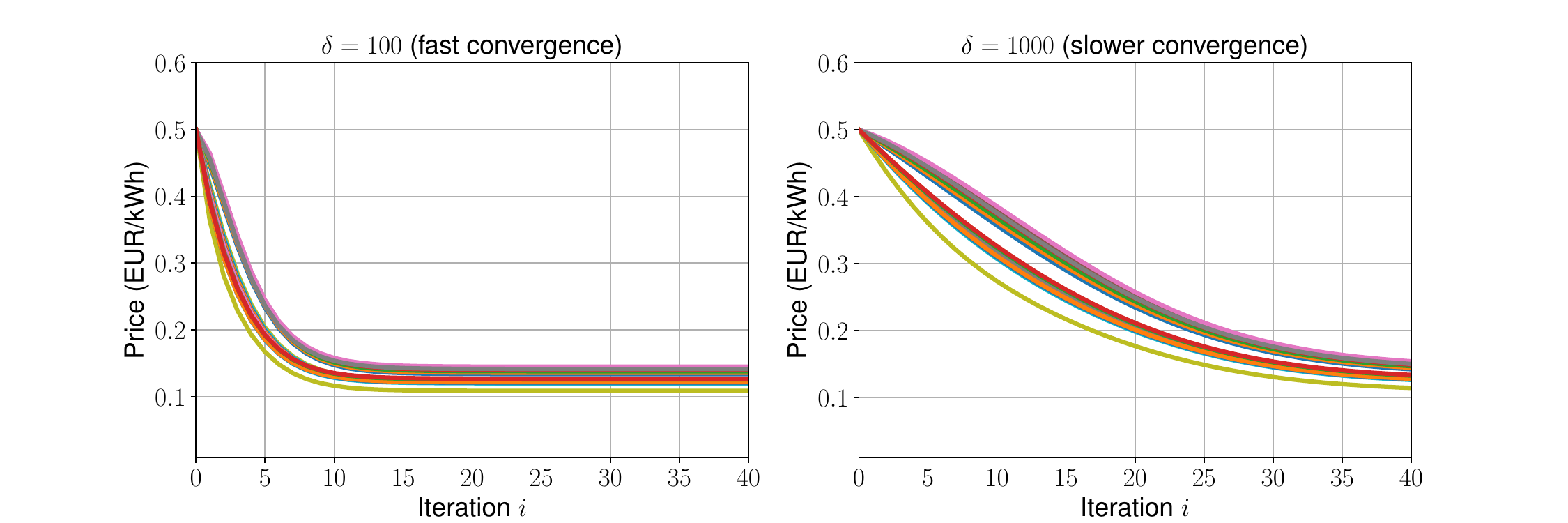}
\caption{Distributed algorithm's performance (Theorem 4 requires $\delta\geq0$) using the Dutch pilot data.}
\label{algorithm2}
\end{figure*}
Figure \ref{Influence_of_T} shows the influence of varying the number of periods on prices, power allocated, revenues, and consumer utilities. We observe that as $T$ increases, the power allocated at each period gets progressively smaller. On the other hand, prices can increase or decrease, depending on the company, and they converge to positive constants. Furthermore, revenues might also increase or decrease, depending on the company (note that the company that achieves the highest revenue is the one that offers the lowest prices, and vice-versa). In view of (\ref{eq1}), the sum of revenues at equilibrium is a constant that matches the sum of all consumer budgets. And hence, whenever the revenue increases (decreases) for a company $k$, at least one other company will incur a loss (gain) in terms of revenue. None of the companies can do better by altering its power availabilities across the time horizon, nor by changing its prices. This follows from the definition of Nash equilibrium. Furthermore, we note that the revenues are proportional to the total capacity, and the company with the highest (lowest) portion of the market is the one that incurs the largest increase (decrease) in revenue. 

Interestingly, in Figure \ref{Influence_of_T} we observe that as $T$ increases, the utilities for consumers also increase, and hence they will be more attracted to demand response programs, which is desirable \cite{DOECOM}. In comparison with the single-period setup \cite{sabita,sabita2}, this shows that the multi-period demand response provides improvements on the consumers' end. This increase, however, does not change significantly beyond a certain number of periods.
To demonstrate the performance of our algorithm, we take the case when $T=1$ and study the algorithm's performance for different values of $\delta$ in Figure \ref{algorithm1}. When $\delta=1000$, we observe that the algorithm converges very fast to the optimal prices and takes about less than $5$ iterations to reach equilibrium. The values are consistent with the values in Figure  \ref{Influence_of_T} when $T=1$, where we used the analytical expressions of the prices. Next, we increase $\delta$ to $10000$ and observe that the algorithm converges at a lower rate, but still fast. Thus, the rate of convergence is inversely proportional to the value of $\delta$. However, when $\delta$ decreases to a negative value, there are no guarantees on convergence. Theorem 4 only guarantees the convergence of the algorithm when $\delta\geq0$. We have verified that our distributed algorithm converges very fast for various values of $\delta$ and alternative values of $T$ and $K$, and the reader might experiment with varying them using our open-source code in \cite{tool2}.

\subsection{Dutch Smart Grid Pilot} {\color{black} To further validate our multi-period-multi-commpany framework, we use data from the Dutch Smart Grid Pilot \cite{dutch}, which was conducted in Zwolle, the Netherlands, for about one year (May, 2014 to May, 2015). Tariffs were announced to consumers a day ahead, and the average consumer behavior was reported}. For a group of $77$ homogeneous consumers, we study the average consumer's demand and payments using experimental prices and the prices derived using our method. Here, we take $K=1$, which is consistent with the Dutch pilot. Also, the experimental prices are in EUR/kWh. 
 
  \subsubsection{Finding the necessary parameters}{\color{black} We find the fixed parameters similarly to the EcoGrid EU experiment. For each consumer $n$, we have {\color{black}$E_n^{\min} \approx 8.8 \ \text{kWh}$.} Then, by (\ref{budget}), we find the minimum necessary daily budget, which is $B_n \approx 1.1\ \text{EUR}$ for each consumer}.

   \subsubsection{Numerical Results} Using the above parameters, we again use (\ref{xx}) and (\ref{p}) to find optimal demands and prices. In Figure \ref{Netherlands_Price_Power}, we plot the average consumer's hourly demand, the prices found experimentally and using the Stackelberg game, and the corresponding total payments by the average consumer. We again observe that our approach leads to smaller prices with a significantly smaller variance. For the average consumer, we observe that significant savings can be achieved using our approach (more than $30\%$). Next, we study the performance of our distributed algorithm in Figure \ref{algorithm2}. As in the case of the EcoGrid EU experimental data, our algorithm achieves fast convergence to optimal prices using only local information.

\section{Generalizations} \label{generalizations}
 
 In the previous sections, we have analyzed our multi-period-multi-company framework under some assumptions to keep the analysis tractable and to reveal various insights on what happens at the equilibrium strategies. Due to the desirable mathematical properties of our framework, it is possible to extend our model at both the consumer-level and company-level. Here, we discuss some such possible extensions. 
 
\subsection{Consumer-Side}
In the utility function (\ref{consumer}), the parameters $\gamma_n$ and $\zeta_n$ for consumer $n$ are time and company independent. However, it is possible, to consider both time-specific and company-specific preferences $\gamma_{n,k,t}$ and $\zeta_{n,k,t}$, which allows consumers to have further flexibilities without violating existence and uniqueness of optimal strategies. In this case, the utility of consumer $n$ would be defined as 
\begin{equation}{u_n(\mathbf{d}_{n})=\sum_{k\in \scr{K}}\sum_{t\in \scr{T}}\gamma_{n,k,t}\ln(\zeta_{n,k,t}+d_{n,k}(t))}. \label{consumer2}\end{equation}
By an analogous analysis to the derivation of (\ref{xx}), it follows that optimal demands, under Assumption \ref{assumption1}, are given by
 \begin{align}
d^*_{n,k}(t)&= \frac{B_n+\sum_{j\in \scr{K}}\sum_{h\in \scr{T}}p_j(h)\zeta_{n,j,h}}{p_k(t)} \Gamma_{n,k,t} \nonumber\\
 &\qquad \qquad\qquad-\zeta_{n,k,t},  \; \forall \;t\in \scr{T}, \; k\in \scr{K},  \label{xx2}
\end{align}
where $$\Gamma_{n,k,t}=\frac{\gamma_{n,k,t}}{\sum_{j\in \scr{K}}\sum_{h\in \scr{T}}\gamma_{n,j,h}}.$$ 
We remark that $\sum_{k\in\mathcal{K}}\sum_{t\in\mathcal{T}}\Gamma_{n,k,t}=1$. Thus, if consumer $n$ prefers a higher demand from company $k$ at time $t$, choosing a higher weight $\gamma_{n,k,t}$ can achieve this. We also note that in (\ref{xx}), where consumer $n$ has identical time-specific/company-specific parameters,  $$\Gamma_{n,k,t}=\frac{1}{KT}.$$ Another alternative, is to expand the constraint set of the optimization problem for consumers to include additional time-specific or company-specific constraints. In general, companies, as leaders of the Stackelberg game, would need to anticipate how consumers would respond to their prices, and given that anticipation, they choose their prices accordingly. Furthermore,  if non-logarithmic utility functions are used by consumers, it might be more difficult to compute a Nash equilibrium for the price-selection game for companies, but the existence of a pure-strategy equilibrium is guaranteed as long as the function 
$$\sum_{t\in \scr{T}}p_k(t)\sum_{n\in \scr{N}}d_{n,k}({\bf{p}}_k,{\bf{p_{-k}}},t)$$
is concave in each $p_k(t)$ over a compact and convex set \cite{basar}, for each company. In case (\ref{consumer2}) is used, by (\ref{xx2}), this condition is satisfied.  

\subsection{Company-Side}

{\color{black} In the current formulation, consumers' demands are coupled through the companies' problems and the power availability constraint (\ref{prob2}). The upper-bound in (\ref{prob2}) is taken to be fixed in the Stackelberg game, but they can be strategically chosen by the power allocation game discussed in Section \ref{game2}. However, this game was solved under restrictive assumptions, such as the absence of network constraints, the full-controllability of generation sources, and the absence of ramping considerations. It is of interest to generalize the power allocation game to alleviate these limitations. Specifically, suppose that power availabilities ${\bf {G}} \in \mathcal{M} \subset \mathbb{R}^{KT}$, where $\mathcal{M}$ represents the transmission and distribution network constraints. One possibility is to assume that $\mathcal{M}$ is a system of linear equations  that approximate power flow equations \cite{eugene,eugene2,eugene3,baran,bolognani}. Furthermore, for simplicity, suppose that company $k$ has a ramping limit $l_{k,t}$ at period $t$. Also, to encode controllability, suppose that $$G^{\min}_{k,t} \leq G_k(t) \leq G^{\max}_{k,t},$$
where $G^{\min}_{k,t}$ ($G^{\max}_{k,t}$) is the minimum (maximum) possible generation  for company $k$ at period $t$.  Thus, company $k$ solves the following optimization problem:
 \begin{eqnarray}
\underset{\mathbf{G}_k }{\hbox{maximize}} && \pi_k({\bf{G}}_k,{\bf{G_{-k}}})
\nonumber \\
\hbox{subject to} && \sum_{t\in \scr{T}}G_k(t) \leq G^{{\rm total}}_k,  \nonumber\\
&& \vert G_k(t)-G_k(t-1)\vert \leq l_{k,t}, \forall t, t-1 \in\scr{T}, \nonumber \\
&& ({\bf{G}}_k,{\bf{G_{-k}}}) \in \mathcal{M}, \label{game3} \\
&& G^{\min}_{k,t} \leq G_k(t) \leq G^{\max}_{k,t},  \ \forall t\in\scr{T}, \nonumber\\
&& G_k(t)\geq0, \ \forall t\in\scr{T}. \nonumber\end{eqnarray}

We have the following result, whose proof is given in the Appendix. 
\begin{theorem}
If the power allocation game (\ref{game3}) is feasible, then, it admits a pure-strategy Nash equilibrium $({\bf{G}}^*_k,{\bf{G^*_{-k}}})$. Furthermore, if $({\bf{G}}^*_k,{\bf{G^*_{-k}}})$ is used for the Stackelberg equilibrium demands and prices given by Theorem \ref{mainTHM}, then, $$\sum_{n\in \scr{N}} d^*_{n,k}(t) = G^*_k(t),\;\; \forall \; t \in \scr{T},\;\; \forall \; k \in \scr{K}.$$ \label{ThmGame3}
\end{theorem}

The above theorem follows from the strict concavity of $\pi_k({\bf{G}}_k,{\bf{G_{-k}}})$ and the compactness and convexity of the constraint set, in addition to the results in Section \ref{game1}. Furthermore, it also demonstrates that it is possible to incentivize consumers to further shift their consumption in a way that is consistent with network considerations and company requirements. Finally, we remark that the control of consumers' demands here is indirect, that is, it is done via the unique equilibrium prices (\ref{p}), which are also affected by consumers' preferences and choices. Hence, at equilibrium, optimal supply provided by companies, ${\bf G^*}$, is equal to aggregate optimal demand, while taking into account consumer budgets and energy needs, in addition to network considerations and company-specific constraints and revenues.

\section{Conclusion and Research Directions} \label{conclusion}
In this paper, we model and solve a novel multi-period-multi-company demand response framework. We formulate a Stackelberg game to capture the interactions between companies and energy consumers, and within the framework of this model, we have obtained optimal prices and demands. Using the closed-form expressions, a power allocation game for companies has been formulated and solved. Furthermore, a distributed algorithm has been proposed to compute all equilibrium strategies using only local information. In the large population regime, an appropriate company-to-user ratio has been derived to maximize the revenue of each individual company.  The paper has shown theoretically and numerically that the multi-period scheme provides more incentives for the participation of energy consumers in demand response management, which is of critical importance \cite{DOECOM}. {\color{black} We have derived a minimum budget condition for consumers that can be used to measure whether or not they are spending more than what is necessary, and case studies using real data reveals potential savings for consumers that can exceed $30\%$}. Numerical studies also demonstrate fast convergence of the proposed distributed algorithm.

While the proposed method focuses on the interplay between competing companies and their consumers, its useful mathematical properties make it generalizable to more consumer-specific and/or company-specific considerations.  For example, it is possible to include period-specific constraints for consumers. The game studied in this paper is multi-period but static. {\color{black} Therefore, it is a one-shot game and all the information are given at the the beginning of the game. Extending it to dynamic information structures, and using tools from dynamic game theory, such as feedback Stackelberg games \cite{basar}, where companies at each period change their prices for the next periods based on the information available at that particular period in which they are making the decisions, is another possible direction. Finally, for the distributed algorithm, it would be interesting to study privacy aspects other than convergence using only local information, such as, the ability of companies to approximate private parameters.}

\section{Acknowledgments}
K. Alshehri thanks King Fahd University of Petroleum and Minerals (KFUPM) for the financial support. Research supported in part by the U.S. Air Force Office of Scientific Research (AFOSR) MURI Grant FA9550-10-1-0573, and in part by the AFOSR Grant FA9550-19-1-0353.

\section{Appendix}
\subsection{Proof of Theorem 1} Note that \begin{equation*}
    B_n \geq \frac{E_n^{\min}+\zeta_nKT}{\sum_{k\in \scr{K}}\sum_{t\in \scr{T}}\frac{1}{KTp_k(t)}}-\zeta_n \sum_{k\in \scr{K}}\sum_{t\in \scr{T}}p_k(t)
\end{equation*}
is the same as \begin{equation*}
\sum_{k\in \scr{K}}\sum_{t\in \scr{T}}\frac{B_n+\zeta_n \sum_{k\in \scr{K}}\sum_{t\in \scr{T}}p_k(t)}{KTp_k(t)}- \sum_{k\in \scr{K}}\sum_{t\in \scr{T}}\zeta_n \geq E_n^{\min}.\end{equation*}
By (\ref{xx}), this is equivalent to 
$\sum_{k\in \scr{K}}\sum_{t\in \scr{T}}d^*_{n,k}(t)\geq \,E^{\min}_n.$

\subsection{Proof of Lemma 1} The matrix $A$ can be represented as 
{\footnotesize{\begin{equation*}\begin{split} \begin{pmatrix}
KT(G_1(1)+Z)& 0&\dots& 0\\
0 & KT(G_1(2)+Z)& \dots& 0\\
\vdots & \ddots\\
0 &\dots&0 & KT(G_K(T)+Z)
\end{pmatrix}\\+\begin{pmatrix}
-Z\\
-Z\\
\vdots \\
-Z
\end{pmatrix}\begin{pmatrix}
1 \dots 1
\end{pmatrix}:=\hat{A}+uv^T\end{split}\end{equation*}}}
Note that $\hat{A}$ is invertible. Furthermore, 
\begin{equation*}
1+v^T\hat{A}^{-1}u=1-\frac{1}{KT}\sum_{k\in \scr{K}}\sum_{t\in \scr{T}}\frac{Z}{G_k(t)+Z}.
\end{equation*}
Since $G_k(t)>0$ and $Z>0$, each element in the summation is less than $1$ and overall value of the summation is less than $KT$, and this clearly leads to $1+v^T\hat{A}^{-1}u\neq0$. By Sherman-Morrison Formula \cite{SM},  if $1+v^T\hat{A}^{-1}u\neq0$, then \begin{equation}A^{-1}=(\hat{A}+uv^T)^{-1}=\hat{A}^{-1}-\frac{\hat{A}^{-1}uv^T\hat{A}^{-1}}{1+v^T\hat{A}^{-1}u}. \label{sm} \end{equation}
Thus, $A$ is invertible and we can apply (\ref{sm}).

\subsection{Proof of Lemma 2} By Lemma 1, the prices are uniquely given by $P=A^{-1}Y$, and by using (\ref{sm}), the price selection for each $k$ at $t$ is 
$$p^*_k(t)=\frac{B}{G_k(t)+Z}\left(\frac{1}{KT-\sum_{j\in \scr{K}}\sum_{h\in \scr{T}}\frac{Z}{G_j(h)+Z}}\right).$$
 Strict positivity follows from $$\frac{B}{G_k(t)+Z}>0  \  \ \text{and} \ \ \sum_{j\in \scr{K}}\sum_{h\in \scr{T}}\frac{Z}{G_j(h)+Z}<KT.$$

\subsection{Proof of Theorem 2}
\begin{itemize}
\item[(i)] By plugging-in the demands given by (\ref{xx}) in the revenue function (\ref{UC}) for $k$, we have 
$$\pi_k=B/K+(Z/K)\sum_{k\in \scr{K}}\sum_{t\in \scr{T}}p_k(t)-Z\sum_{t\in \scr{T}}p_k(t),$$

\noindent
which is concave (linear) in each $p_k(t)$. Thus, by the compactness of $\scr{L}_{k,t}$,  there exists a pure-strategy Nash Equilibrium (NE) \cite{basar}. Next, suppose that a company $k$ deviates from (\ref{p}) and announces a price of $\hat{p}_k(t)=p^*_k(t)+\epsilon$ at a fixed time $t$. If $\epsilon>0$, then
    $$\hat{\pi}-\pi_k=\epsilon  \frac{Z-ZK}{K} \leq 0,$$
    where the inequality follows from $ZK\geq K$. Thus, $k$ has no incentive to increase the prices from (\ref{p}). Furthermore, since the prices given by (\ref{p}) are attained the equality of the capacity constraint in (\ref{prob2}), company $k$ has no incentive to choose $\epsilon <0$ because it will not result in selling more energy. Therefore, for every period $t$, company $k$ does not benefit from deviating from (\ref{p}). Hence, the prices given by (\ref{p}) maximize the revenues and constitute a NE.


\item[(ii)] By the uniqueness of the demands given by (\ref{xx}) and using (i), it follows that there exists a unique Stackelberg equilibrium and it is given by the pair ${\bf{d^*(p)}}$ and (\ref{p}).
\end{itemize}
\subsection{Proof of Theorem 3} \label{game2proof}
Note that the revenue $\pi_k({\bf{G}}_k,{\bf{G_{-k}}})$ is equivalent  to
\begin{equation} \sum_{t\in \scr{T}}\frac{BG_k(t)}{(G_k(t)+N)(\alpha_{-k}-\sum_{h\in \scr{T}}\frac{N}{G_k(h)+N})} ,\label{Uk2} \end{equation}
where   $$\alpha_{-k}:=  KT-\sum_{j\in \scr{K},j\neq k\,\,} \sum_{h\in \scr{T}}\frac{N}{G_j(h)+N}> T.$$ 
Note that $\alpha_{-k}$ depends on the strategies of other companies and it is fixed for company $k$. 
A pure-strategy Nash equilibrium exists if $\pi_k$ is concave in each $G_k(t)\in\scr{P}_{k,t}$ for each company $k$ and if  $\scr{P}_{k,t}$ is a compact subset of $\mathbb{R}$ \cite{basar}. 
Since it is clear that $\scr{P}_{k,t}$ is compact, it is enough to show concavity of ${ \rm company,k}$. From (\ref{Uk2}), via a sequence of mathematical tricks, 
\begin{align}\pi_k&=\frac{B\sum_{t\in \scr{T}}G_k(t)\prod_{h\neq t}(G_k(h)+N)\frac{G_k(t)+N}{G_k(t)+N}}{\prod_{h\in\scr{T}}(G_k(h)+N)(\alpha_{-k}-\sum_{h\in \scr{T}}\frac{N}{G_k(h)+N})} \nonumber \\
&=\frac{B\prod_{h\in\scr{T}}(G_k(h)+N)\sum_{t\in \scr{T}}\frac{G_k(t)}{G_k(t)+N}}{\prod_{h\in\scr{T}}(G_k(h)+N)(\alpha_{-k}-\sum_{h\in \scr{T}}\frac{N}{G_k(h)+N})} \nonumber \\
&=\frac{B\sum_{t\in \scr{T}}\frac{G_k(t)}{G_k(t)+N}}{\alpha_{-k}-\sum_{h\in \scr{T}}\frac{G_k(h)+N}{G_k(h)+N}+\sum_{t\in \scr{T}}\frac{G_k(t)}{G_k(t)+N}} \nonumber\\
&=B\frac{\sum_{t\in \scr{T}}\frac{G_k(t)}{G_k(t)+N}}{(\alpha_{-k}-T)+\sum_{t\in \scr{T}}\frac{G_k(t)}{G_k(t)+N}} \nonumber \\ 
&=:\frac{f}{\gamma_{-k}+f}.\end{align}
Note that $f_{G_k(t)}=\frac{\partial f}{\partial G_k(t)}=\frac{N}{(G_k(t)+N)^2} > 0$ and
\begin{equation*}\frac{\partial \pi_k}{\partial G_k(t)}=\frac{f_{G_k(t)}\gamma_{-k}}{(\gamma_{-k}+f)^2}=\frac{N\gamma_{-k}}{(\gamma_{-k}+f)^2(G_k(t)+N)^2} >0. \label{der1} \end{equation*}
This leads to
\begin{equation*}\frac{\partial^2 \pi_k}{\partial G_k(t)^2}=\frac{[-2N\gamma_{-k}][(\gamma_{-k}+f)+f_{G_k(t)}(G_k(t)+N)]}{[(\gamma_{-k}+f)(G_k(t)+N)]^2}, \label{der2} \end{equation*}
which is strictly negative since $f, f_{G_k(t)} ,N,\gamma_{-k} >0$.
Hence, strict concavity holds. We can relax the non-negativity constraint as the solution will be positive by the properties of the objective function. The Lagrange function for company $k$ is then given by
\begin{equation}
L_{k}({\bf{G}}_k,{\bf{G_{-k}}},\lambda_k)= \pi_k + \lambda_k \left(\sum_{t\in \scr{T}}G_k(t)-G^{{\rm total}}_k\right),
\end{equation}
and by the first-order necessary condition $\nabla L=0$,
\begin{eqnarray}\lambda_k &=&-\frac{N\gamma_{-k}}{(\gamma_{-k}+f)^2(G_k(t)+N)^2},\,\,\,  \forall \,t\in \scr{T}\\
\frac{\partial L_k}{\partial \lambda_k} &=& 0 \implies \sum_{t\in \scr{T}}G_k(t)=G^{{\rm total}}_k.
\end{eqnarray}
Thus, for company $k$, elements of ${\bf G}_k$ must be identical, and must add up to $G^{{\rm total}}_k$.
\subsection{Proof of Theorem 4}
 To find an appropriate $\epsilon^{(i)}_{k,t}$ that leads to the convergence, recall that the prices must be positive. The algorithm diverges whenever any $p^{(i)}_k(t)$ is negative, which might happen when  $\sum_{n\in \scr{N}}d^{(i)}_{n,k}(t)<G_k(t)$, for any company $k\in\scr{K}$ at any time $t\in\scr{T}$ in iteration $i$. To avoid this, it suffices to require 
$ p^{(i)}_k(t)\epsilon^{(i)}_{k,t} > \left |  \sum_{n\in \scr{N}}d^{(i)}_{n,k}(t)-G_k(t)  \right |$  
whenever we have $\sum_{n\in \scr{N}}d^{(i)}_{n,k}(t)<G_k(t)$. This translates into requiring $$p^{(i)}_k(t)\epsilon^{(i)}_{k,t}> G_k(t)-\sum_{n\in \scr{N}}d^{(i)}_{n,k}(t) $$ for any $k\in\scr{K}$, $t\in\scr{T}$, and~$i$. By (\ref{xx}), it follows that we need
\begin{equation}
\epsilon^{(i)}_{k,t}
>\frac{G_k(t)-\sum_{n\in \scr{N}}\left(\frac{B_n+\sum_{j\in \scr{K}}\sum_{h\in \scr{T}}p^{(i)}_j(h)}{KTp^{(i)}_k(t)}-1\right)}{p^{(i)}_k(t)}.
\label{bound}\end{equation}
The bound (\ref{bound}) is the tightest one, but using it to find $\epsilon^{(i)}_{k,t}$ is not implementable. By choosing 
 \begin{equation}\epsilon^{(i)}_{k,t} \geq \frac{G_k(t)+N}{p^{(i)}_k(t)},\label{eps}\end{equation}
 convergence is guaranteed since $$\frac{B_n+\sum_{j\in \scr{K}}\sum_{h\in \scr{T}}p^{(i)}_j(h)}{KTp^{(i)}_k(t)}\geq0.$$

\subsection{Proof of Theorem 5} Suppose that $$\frac{K}{N}<\frac{B_n}{p^{{\rm max}}TG^*_k(t)}.$$
By (\ref{pppp}), this implies that $$p^{{\rm max}}<\frac{NB_n}{KTG^*_k(t)}=p^*_k(t),\,\,\, \forall \,t\in \scr{T},\,\, \forall \,k\in \scr{K},$$
and companies will charge $p^{{\rm max}}$, which implies  $$\sum_{k\in \scr{K}}\pi_k=p^{{\rm max}}KTG^*_k(t)<NB_n=\sum_{n\in \scr{N}}B_n,$$
which means that the sum of the revenues is strictly less than the sum of the budgets and hence companies incur losses, compared to the  equilibrium prices. On the other hand, $$\frac{K}{N} \geq \frac{B_n}{p^{{\rm max}}TG^*_k(t)} $$ is equivalent to  $$p^{{\rm max}}\geq \frac{NB_n}{KTG^*_k(t)}=p^*_k(t),\,\,\, \forall \,t\in \scr{T},\,\, \forall \,k\in \scr{K}.$$ Furthermore, we have $$\sum_{k\in \scr{K}}\pi_k=p^*_k(t)KTG^*_k(t)=NB_n=\sum_{n\in \scr{N}}B_n.$$

{\color{black} 
\subsection{Proof of Theorem 6}
From Section \ref{game2proof}, the revenue function $$\pi_k({\bf{G}}_k,{\bf{G_{-k}}})$$ is strictly concave in $G_k(t)$ for each company $k$ at period $t$. Furthermore, since the constraint set is convex and compact in $G_k(t)$, existence of a pure-strategy Nash equilibrium is guaranteed \cite{basar}. The rest of the proof readily follows from Theorem 2.
}

\bibliographystyle{IEEEtran}
\bibliography{references}

\begin{thebibliography}{10}
\providecommand{\url}[1]{#1}
\csname url@samestyle\endcsname
\providecommand{\newblock}{\relax}
\providecommand{\bibinfo}[2]{#2}
\providecommand{\BIBentrySTDinterwordspacing}{\spaceskip=0pt\relax}
\providecommand{\BIBentryALTinterwordstretchfactor}{4}
\providecommand{\BIBentryALTinterwordspacing}{\spaceskip=\fontdimen2\font plus
\BIBentryALTinterwordstretchfactor\fontdimen3\font minus
  \fontdimen4\font\relax}
\providecommand{\BIBforeignlanguage}[2]{{%
\expandafter\ifx\csname l@#1\endcsname\relax
\typeout{** WARNING: IEEEtran.bst: No hyphenation pattern has been}%
\typeout{** loaded for the language `#1'. Using the pattern for}%
\typeout{** the default language instead.}%
\else
\language=\csname l@#1\endcsname
\fi
#2}}
\providecommand{\BIBdecl}{\relax}
\BIBdecl

\bibitem{sabita}
S.~Maharjan, Q.~Zhu, Y.~Zhang, S.~Gjessing, and T.~Ba\c{s}ar, ``Dependable
  demand response management in the smart grid: a {S}tackelberg game
  approach,'' \emph{IEEE Trans. Smart Grid}, vol.~4, no.~1, pp. 120--132, 2013.

\bibitem{DR}
M.~H. Albadi and E.~F. El-Saadany, ``Demand response in electricity markets: An
  overview,'' \emph{IEEE Power Engineering Society General Meeting}, 2007.

\bibitem{DSM}
P.~Palensky and D.~Dietrich, ``Demand side management: Demand response,
  intelligent energy systems, and smart loads,'' \emph{IEEE Trans. Industrial
  Informatics}, vol.~7, no.~3, pp. 381--388, 2011.

\bibitem{DSview}
D.~S. Kirschen, ``Demand-side view of electricity markets,'' \emph{IEEE Trans.
  Power Systems}, vol.~18, no.~2, pp. 520--527, 2003.

\bibitem{DReff}
K.~Spees and L.~B. Lave, ``Demand response and electricity market efficiency,''
  \emph{The Electricity Journal}, vol.~20, no.~3, pp. 69--85, 2007.

\bibitem{challengesDR}
S.~Nolan and M.~O'Malley, ``Challenges and barriers to demand response
  deployment and evaluation,'' \emph{Applied Energy}, vol. 152, pp. 1--10,
  2015.

\bibitem{DRprice}
S.~Chan, K.~Tsui, H.~Wu, Y.~Hou, and F.~F. Wu, ``Load/price forecasting and
  managing demand response for smart grids: methodologies and challenges,''
  \emph{IEEE Signal Processing Magazine}, vol.~29, no.~5, pp. 68--85, 2012.

\bibitem{surveyDRM}
J.~S. Vardakas, N.~Zorba, and C.~V. Verikoukis, ``A survey on deamnd response
  programs in smart grids: pricing methods and optimization algorithms,''
  \emph{IEEE Communications Surveys and Tutorials}, vol.~17, no.~1, pp.
  152--178, 2013.

\bibitem{reviewAE}
J.~Wang, H.~Zhong, Z.~Ma, Q.~Xia, and C.~Kang, ``Review and prospect of
  integrated demand response in the multi-energy system,'' \emph{Applied
  Energy}, vol. 202, pp. 772--782, 2017.

\bibitem{loadadapt}
P.~B. Luh, Y.-C. Ho, and R.~Muralidharan, ``Load adaptive pricing: an emerging
  tool for electric utilities,'' \emph{IEEE Trans. Automatic Control}, vol.~27,
  no.~2, pp. 320--329, 1982.

\bibitem{survey}
W.~Saad, Z.~Han, H.~V. Poor, and T.~Ba\c{s}ar, ``Game-theoretic methods for the
  smart grid: an overview of microgrid systems, demand-side management, and
  smart grid communications,'' \emph{IEEE Signal Processing Magazine}, vol.~29,
  no.~5, pp. 86--105, Sept. 2012.

\bibitem{walidPHEV}
------, ``A noncooperative game for double auction-based energy trading between
  {PHEV}s and distribution grids,'' \emph{2011 IEEE International Conference on
  Smart Grid Communications (SmartGridComm)}, pp. 267--272, 2011.

\bibitem{trading}
Y.~Wang, W.~Saad, Z.~Han, H.~V. Poor, and T.~Ba\c{s}ar, ``A game-theoretic
  approach to energy trading in the smart grid,'' \emph{IEEE Trans. Smart
  Grid}, vol.~5, no.~3, pp. 1439--1450, 2014.

\bibitem{walidPHEV2}
W.~Tushar, W.~Saad, H.~V. Poor, and D.~B. Smith, ``Economics of electric
  vehicle charging: A game theoretic approach,'' \emph{IEEE Trans. Smart Grid},
  vol.~3, no.~4, pp. 1767--1778, 2012.

\bibitem{gao}
B.~Gao, W.~Zhang, Y.~Tang, M.~Hu, M.~Zhu, and H.~Zhan, ``Game-theoretic energy
  management for residential users with dischargeable plug-in electric
  vehicles,'' \emph{Energies}, vol.~7, no.~11, pp. 7499--7518, 2014.

\bibitem{twolevel}
B.~Chai, J.~Chen, Z.~Yang, and Y.~Zhang, ``Demand response management with
  multiple utility companies: A two-level game approach,'' \emph{IEEE Trans.
  Smart Grid}, vol.~5, no.~2, pp. 722--731, 2014.

\bibitem{tansu}
E.~Nekouei, T.~Alpcan, and D.~Chattopadhyay, ``A game-theoretic analysis of
  demand response in electricity markets,'' \emph{2014 IEEE PES General
  Meeting}, 2014.

\bibitem{yaagoubi}
N.~Yaagoubi and H.~T. Mouftah, ``A distributed game theoretic approach to
  energy trading in the smart grid,'' \emph{2015 IEEE Electrical Power and
  Energy Conference (EPEC)}, pp. 203--208, 2015.

\bibitem{tushar2}
W.~Tushar, B.~Chai, C.~Yuen, D.~B. Smith, K.~L. Wood, Z.~Yang, and H.~V. Poor,
  ``Three-party energy management with distributed energy resources in smart
  grid,'' \emph{IEEE Trans. Industrial Informatics}, vol.~62, no.~4, 2015.

\bibitem{sabita2}
S.~Maharjan, Q.~Zhu, Y.~Zhang, S.~Gjessing, and T.~Ba\c{s}ar, ``Demand response
  management in the smart grid in a large population regime,'' \emph{IEEE
  Trans. Smart Grid}, vol.~7, no.~1, pp. 189--199, 2016.

\bibitem{sabita3}
S.~Maharjan, Y.~Zhang, S.~Gjessing, and D.~H. Tsang, ``User-centric demand
  response management in the smart grid with multiple providers,'' \emph{IEEE
  Trans. Emerging Topics in Computing}, vol.~4, no.~5, 2014.

\bibitem{han}
K.~Han, J.~Lee, and J.~Choi, ``Evaluation of demand-side management over
  pricing competition of multiple suppliers having heterogeneous energy
  sources,'' \emph{Energies}, vol.~10, no.~9, 2017.

\bibitem{amir}
A.-H. Mohsenian-Rad, V.~W.~S. Wong, J.~Jatskevich, R.~Schober, and
  A.~Leon-Garcia, ``Autonomous demand-side management based on game-theoretic
  energy consumption scheduling for the future smart grid,'' \emph{IEEE Trans.
  Smart Grid}, vol.~1, no.~3, pp. 320--331, 2010.

\bibitem{hazem}
H.~Soliman and A.~Leon-Garcia, ``Game-theoretic demand-side management with
  storage devices for the future smart grid,'' \emph{IEEE Trans. Smart Grid},
  vol.~5, no.~3, pp. 1475--1485, 2014.

\bibitem{roh}
H.-T. Roh and J.-W. Lee, ``Residential demand response scheduling with
  multiclass appliances in the smart grid,'' \emph{IEEE Trans. Smart Grid},
  vol.~7, no.~1, 2016.

\bibitem{zhudiff}
Q.~Zhu, Z.~Han, and T.~Ba\c{s}ar, ``A differential game approach to distributed
  demand side management in smart grid,'' \emph{2012 IEEE International
  Conference on Communications}, pp. 3345--3350, 2012.

\bibitem{PAR}
H.~K. Nguyen, J.~B. Song, and Z.~Han, ``Demand side management to reduce
  {p}eak-to-{a}verage ratio using game theory in smart grid,'' \emph{1st IEEE
  INFOCOM Workshop on Communications and Control for Sustainable Energy
  Systems: Green Networking and Smart Grids}, pp. 91--96, 2012.

\bibitem{collins}
L.~D. Collins and R.~H. Middleton, ``Distributed demand peak reduction with
  non-cooperative players and minimal communication,'' \emph{IEEE Trans. Smart
  Grid}, vol.~PP, no.~99, 2017.

\bibitem{repeated}
L.~Song, Y.~Xiao, and M.~van~der Schaar, ``Demand side management in smart
  grids using a repeated game framework,'' \emph{IEEE Journal of Selected Areas
  in Communications}, vol.~32, no.~7, pp. 1412--1424, 2014.

\bibitem{fourstage}
S.~Bu and F.~R. Yu, ``A game-theoretical scheme in the smart grid with
  demand-side management: towards a smart cyber-physical power
  infrastructure,'' \emph{IEEE Trans. Emerging Topics in Computing}, vol.~1,
  no.~1, pp. 22--32, 2013.

\bibitem{dayahead}
L.~Jia and L.~Tong, ``Day ahead dynamic pricing for demand response in dynamic
  environments,'' \emph{52nd IEEE Conference on Decision and Control}, pp.
  5608--5613, 2013.

\bibitem{wei}
W.~Wei, F.~Liu, and S.~Mei, ``Energy pricing and dispatch for smart grid
  retailers under demand response and market price uncertainty,'' \emph{IEEE
  Trans. Smart Grid}, vol.~6, no.~3, 2015.

\bibitem{mywork}
K.~Alshehri, J.~Liu, X.~Chen, and T.~Ba\c{s}ar, ``A {S}tackelberg game for
  multi-period demand response management in the smart grid,''
  \emph{Proceedings of the 54th IEEE Conference on Decision and Control}, pp.
  5889--5894, 2015.

\bibitem{basar}
T.~Ba\c{s}ar and G.~J. Olsder, \emph{Dynamic Noncooperative Game Theory}.\hskip
  1em plus 0.5em minus 0.4em\relax SIAM, 1999.

\bibitem{srikant}
R.~Srikant, \emph{The Mathematics of Internet Congestion Control}.\hskip 1em
  plus 0.5em minus 0.4em\relax Birkh{\"a}user, 2004.

\bibitem{shadow}
F.~Kelly, A.~Maulloo, and D.~Tan, ``Rate control for communication networks:
  shadow prices, proportional fairness and stability,'' \emph{Journal of the
  Operational Research Society}, vol.~49, no.~3, pp. 237--252, 1998.

\bibitem{basarDR}
T.~Ba\c{s}ar and R.~Srikant, ``Revenue-maximizing pricing and capacity
  expansion in a many-users regime,'' \emph{IEEE International Conference on
  Computer Communications}, vol.~1, pp. 294--301, 2002.

\bibitem{basarDR2}
------, ``A {S}tackelberg network game with a large number of followers,''
  \emph{Journal of Optimization Theory and Applications}, vol. 115, no.~3, pp.
  479--490, 2002.

\bibitem{fan2}
Z.~Fan, ``A distributed demand response algorithm and its application to {PHEV}
  charging in smart grids,'' \emph{IEEE Trans. Smart Grid}, vol.~3, no.~3, pp.
  1280--1290, 2012.

\bibitem{DRadaptation}
------, ``Distributed demand response and user adaptation in smart grids,''
  \emph{2011 IFIP/IEEE International Symposium on Integrated Network Management
  (IM)}, 2011.

\bibitem{NL}
D.~P. Bertsekas, \emph{Nonlinear Programming}.\hskip 1em plus 0.5em minus
  0.4em\relax Athena Scientific, 2008.

\bibitem{bosebook}
S.~Bose and S.~H. Low, \emph{Smart Grid Control: Overview and Research
  Opportunities}, ser. Power Electronics and Power Systems.\hskip 1em plus
  0.5em minus 0.4em\relax Springer, 2018, ch. Some Emerging Challenges in
  Electricity Markets, pp. 29--45.

\bibitem{dutch}
E.~A.~M. Klaassenab, C.~B.~A. Kobus, J.~Frunt, and J.~G. Slootweg,
  ``Responsiveness of residential electricity demand to dynamic tariffs:
  Experiences from a large field test in the netherlands,'' \emph{Applied
  Energy}, vol. 183, pp. 1065--1074, 2016.

\bibitem{ecogrid}
P.~Lund, P.~Nyeng, R.~D. Grandal, S.~H. S{\o}rensen, M.~F. Bendtsen, G.~le~Ray,
  E.~M. Larsen, J.~Mastop, F.~Judex, F.~Leimgruber, K.~J. Kok, and P.~A.
  MacDougall, ``Eco{G}rid {EU} deliverable 6.7: overall evaluations and
  conclusions,'' Eco{G}rid {EU}, Tech. Rep., 2016.

\bibitem{tool2}
\BIBentryALTinterwordspacing
K.~Alshehri, ``Multi-period demand response,'' \emph{{G}it{H}ub repository}.
  [Online]. Available: \url{https://github.com/kalsheh2/DemandResponse}
\BIBentrySTDinterwordspacing

\bibitem{IFAC}
K.~Alshehri, S.~Bose, and T.~Ba\c{s}ar, ``Cash-settled options for wholesale
  electricity markets,'' \emph{Proc. 20th IFAC World Congress (IFAC WC 2017),
  Toulouse, France}, pp. 14\,147--14\,153, 2017.

\bibitem{DOECOM}
U.~S. DOE, \emph{Communications Requirements of Smart Grid Technologies}, US
  Department of Energy, 2010.

\bibitem{eugene}
D.~Bertsimas, E.~Litvinov, X.~A. Sun, J.~Zhao, and T.~Zheng, ``Adaptive robust
  optimization for the security constrained unit commitment problem,''
  \emph{IEEE Transactions on Power Systems}, vol.~28, no.~1, pp. 52--63, 2013.

\bibitem{eugene2}
D.~Gan and E.~Litvinov, ``Energy and reserve market designs with explicit
  consideration to lost opportunity costs,'' \emph{IEEE Transactions on Power
  Systems}, vol.~18, no.~1, pp. 53--59, 2003.

\bibitem{eugene3}
E.~Litvinov, T.~Zheng, G.~Rosenwald, and P.~Shamsollahi, ``Marginal loss
  modeling in lmp calculation,'' \emph{IEEE Transactions on Power Systems},
  vol.~19, no.~2, pp. 880--888, 2004.

\bibitem{baran}
M.~Baran and F.~Wu, ``Network reconfiguration in distribution systems for loss
  reduction and load balancing,'' \emph{IEEE Transactions on Power Systems},
  vol.~4, no.~2, pp. 1401--1407, 1989.

\bibitem{bolognani}
S.~Bolognani and S.~Zampieri, ``On the existence and linear approximation of
  the power flow solution in power distribution networks,'' \emph{IEEE
  Transactions on Power Systems}, vol.~31, no.~1, pp. 163--172, 2016.

\bibitem{SM}
J.~Sherman and W.~J. Morrison, ``Adjustment of an inverse matrix corresponding
  to a change in one element of a given matrix,'' \emph{Annals of Mathematical
  Statistics}, vol.~21, no.~1, pp. 124--127, 1950.

\end{thebibliography}

\end{document}